\documentclass[preprinst,11pt]{elsarticle}\usepackage{amsmath,amssymb,amsfonts,amscd,mathrsfs}
\usepackage{graphicx,fancybox,color,yhmath}
\usepackage{arydshln,MnSymbol}
\usepackage{pdfsync,upgreek,dsfont}
\usepackage{enumerate}
\usepackage[utf8]{inputenc}
\usepackage[normalem]{ulem}
\newcommand{\ls}[2]{#2}
\newcommand{\nz}{{n_{z}}}
\newcommand{\nw}{{n_{w}}}
\newcommand{\ny}{{n_{y}}}

\newcommand{\sy}{G}
\newcommand{\mi}{M}

\newcommand{\col}{\text{col}}
\newcommand{\pl}{\bul}
\newcommand{\Le}{L}

\newcommand{\De}{\c{D}}
\newcommand{\X}{{X}}

\newcommand{\Hi}{H}

\newcommand{\rhi}{R\Hi_\infty}
\newcommand{\rli}{R\Le_\infty}

\newcommand{\st}{\ |\ }

\newcommand{\Ci}{\C_=^\infty}
\newcommand{\eps}{\epsilon}

\newcommand{\Ag}{A}
\newcommand{\Bg}{B}
\newcommand{\Cg}{C}
\newcommand{\Dg}{D}

\newcommand{\Cgr}{C_F}
\newcommand{\Dgr}{D_F}
\newcommand{\Ddgr}{E}
\newcommand{\Ddp}{D_d}

\newcommand{\Af}{\c A}
\newcommand{\Bf}{\c B}
\newcommand{\Cf}{\c C}
\newcommand{\Df}{\c D}
\newcommand{\Xf}{\c X}
\newcommand{\Zf}{\c Z}
\newcommand{\Y}{Y}

\newcommand{\arc}[1]{\arraycolsep#1ex}

\newcommand{\skipthis}[1]{}
\newcommand{\bm}[1]{ \mbox{\boldmath$ #1 $} }

\renewcommand{\t}[1]{ \tilde{#1} }

\newcommand{\h}[1]{ \widehat{#1} }

\renewcommand{\c}[1]{{\cal #1}}

\newcommand{\epro}{ \hfill\mbox{\rule{2mm}{2mm}} }
\newcommand{\Cl}{ {\mathbb{C}^-} }

\newcommand{\R}{ {\mathbb{R}} }
\newcommand{\N}{ {\mathbb{N}} }

\newcommand{\diag}{ \operatornamewithlimits{diag} }

\renewcommand{\r}[1]{(\ref{#1})}

\newtheorem{theo}{Theorem}

\newtheorem{hypo}[theo]{Hypothesis}
\newtheorem{lemm}[theo]{Lemma}
\newtheorem{coro}[theo]{Corollary}
\newtheorem{defi}[theo]{Definition}
\newtheorem{rema}[theo]{Remark}

\newcommand{\theorem}[1]{\begin{theo}#1\end{theo}}

\newcommand{\lemma}[1]{\begin{lemm}#1\end{lemm}}
\newcommand{\definition}[1]{\begin{defi}#1\end{defi}}

\newcommand{\mun}[1]{\begin{multline*}#1\end{multline*}}

\newcommand{\equ}[1]{\begin{equation}#1\end{equation}}
\newcommand{\eqn}[1]{$$#1$$}
\newcommand{\eql}[2]{\begin{equation}\label{#1}#2\end{equation}}

\newcommand{\ear}[1]{\begin{eqnarray}#1\end{eqnarr{{{}}}ay}}

\newcommand{\arr}[2]{\begin{array}{#1}#2\end{array}}
\newcommand{\mas}[2]{\left[\begin{array}{#1}#2\end{array}\right]}
\newcommand{\mat}[2]{\left(\begin{array}{#1}#2\end{array}\right)}

\newcommand{\ga}{\gamma}
\newcommand{\om}{\omega}
\newcommand{\io}{i\omega}
\newcommand{\te}[1]{\text{\ \ #1\ \ }}
\newcommand{\kyp}[3]{\mathscr{L}\left(#1,#2,#3\right)}
\newcommand{\be}{\beta}

\newcommand{\al}{\alpha}

\newcommand{\cl}{\prec}
\newcommand{\cg}{\succ}
\newcommand{\cge}{\succeq}
\newcommand{\cle}{\preceq}

\catcode`\=13
\def{\relax}
\newcommand{\hl}{\\\hline}

\renewcommand{\h}{\hat}

\newcommand{\remark}[1]{\begin{rema}#1\end{rema}}

\usepackage{tikz,pgfplots,datatool}
\usetikzlibrary{positioning,calc,arrows,shapes,shadows}
\tikzset{
auto,
sys/.style 2 args={
rectangle,
draw,
rounded corners,
drop shadow,
fill=white,
minimum height=#2,
minimum width=#1,
inner sep=\dn},
sum/.style={circle,draw,draw=black,inner sep=0mm,minimum size=2mm},
jun/.style={circle,draw,draw=black,inner sep=0mm,minimum size=0mm},
>={latex},
every path/.style={rounded corners},
}

\newcommand{\tio}[4]{\coordinate (#1) at ($(#2.south #3)!#4!(#2.north #3)$)}

\def\dn{1ex}
\def\dl{3*\dn}

\tikzstyle{sy0}=[sys={0*\dn}{0*\dn}]
\tikzstyle{sy1}=[sys={12*\dn}{8*\dn}]
\tikzstyle{sy2}=[sys={8*\dn}{6*\dn}]
\tikzstyle{sy3}=[sys={5*\dn}{5*\dn}]

\usepackage{xparse} 

\DeclareDocumentCommand \ucsysp {o m o m o m}{%
\IfNoValueTF{#1}{%
\node[sy1](s) at (0,0) {#2};%
}{%
\expandafter\node\expandafter[#1](s) at (0,0) {#2};%
}%
\IfNoValueTF{#3}{%
\node[sy3, above = 1*\dl of s](d){#4};%
}{%
\expandafter\node\expandafter[#3, above = 1*\dl of s](d){#4};%
}%
\IfNoValueTF{#5}{%
\node[sy3, below = 1*\dl of s](c){#6};%
}{%
\expandafter\node\expandafter[#5, below = 1*\dl of s](c){#6};%
}%
\ucsyspcontinued
}

\DeclareDocumentCommand \ucsyspcontinued {m m m m m m o}{%
\IfNoValueTF{#7}{%
\def\temp{1}%
}{%
\def\temp{#7}%
}%
\tio{i1}{s}{east}{1/4};
\tio{i2}{s}{east}{2/4};
\tio{i3}{s}{east}{3/4};

\tio{o1}{s}{west}{1/4};
\tio{o2}{s}{west}{2/4};
\tio{o3}{s}{west}{3/4};

\draw[->] (o3) -- ($(o3) - (\temp*\dl, 0)$) |- node[pos = 0.25, left]{#1}(d);
\draw[<-] (i3) -- ($(i3) + (\temp*\dl, 0)$) |- node[pos = 0.25, right]{#2}(d);
\draw[->] (o1) -- ($(o1) - (\temp*\dl, 0)$) |- node[pos = 0.25, left]{#3}(c);
\draw[<-] (i1) -- ($(i1) + (\temp*\dl, 0)$) |- node[pos = 0.25, right]{#4}(c);

\draw[->] (o2) -- node[at end, above]{#5}($(o2) - (\temp*\dl + 1*\dl, 0)$);
\draw[<-] (i2) -- node[at end, above]{#6}($(i2) + (\temp*\dl + 1*\dl, 0)$);
}

\DeclareDocumentCommand \ucsys {o m o m o m}{%
\IfNoValueTF{#1}{%
\node[sy1](s) at (0,0) {#2};%
}{%
\expandafter\node\expandafter[#1](s) at (0,0) {#2};%
}%
\IfNoValueTF{#3}{%
\node[sy3, above = 1*\dl of s](d){#4};%
}{%
\expandafter\node\expandafter[#3, above = 1*\dl of s](d){#4};%
}%
\IfNoValueTF{#5}{%
\node[sy3, below = 1*\dl of s](c){#6};%
}{%
\expandafter\node\expandafter[#5, below = 1*\dl of s](c){#6};%
}%
\ucsyscontinued
}

\DeclareDocumentCommand \ucsyscontinued {m m m m o}{%
\IfNoValueTF{#5}{%
\def\tempa{1}%
}{%
\def\tempa{#5}%
}%
\tio{i1}{s}{east}{1/3};
\tio{i3}{s}{east}{2/3};

\tio{o1}{s}{west}{1/3};
\tio{o3}{s}{west}{2/3};

\draw[->] (o3) -- ($(o3) - (\tempa*\dl, 0)$) |- node[pos = 0.25, left]{#1}(d);
\draw[<-] (i3) -- ($(i3) + (\tempa*\dl, 0)$) |- node[pos = 0.25, right]{#2}(d);

\draw[->] (o1) -- ($(o1) - (\tempa*\dl, 0)$) |- node[pos = 0.25, left]{#3}(c);
\draw[<-] (i1) -- ($(i1) + (\tempa*\dl, 0)$) |- node[pos = 0.25, right]{#4}(c);
}

\DeclareDocumentCommand \usys {o m o m m m m m o}{%
\def\tempa{#1}%
\def\tempb{#3}%
\IfNoValueTF #1 {%
\def\tempa{sy1}%
}{}%
\IfNoValueTF #3 {%
\def\tempb{sy3}%
}{}%
\IfNoValueTF{#9}{%
\def\tempc{1}%
}{%
\def\tempc{#9}%
}%
\expandafter\node\expandafter[\tempa](s) at (0,0) {#2};
\expandafter\node\expandafter[\tempb, above = 1*\dl of s](d){#4};

\tio{i1}{s}{east}{1/3};
\tio{i3}{s}{east}{2/3};

\tio{o1}{s}{west}{1/3};
\tio{o3}{s}{west}{2/3};

\draw[->] (o3) -- ($(o3) - (\tempc*\dl, 0)$) |- node[pos = 0.25, left]{#5}(d);
\draw[<-] (i3) -- ($(i3) + (\tempc*\dl, 0)$) |- node[pos = 0.25, right]{#6}(d);

\draw[->] (o1) -- node[at end, above]{#7}($(o1) - (\tempc*\dl + \dl, 0)$);
\draw[<-] (i1) -- node[at end, above]{#8}($(i1) + (\tempc*\dl + \dl, 0)$);
}

\DeclareDocumentCommand \csys {o m o m m m m m o}{%
\def\tempa{#1}%
\def\tempb{#3}%
\IfNoValueTF #1 {%
\def\tempa{sy1}%
}{}%
\IfNoValueTF #3 {%
\def\tempb{sy3}%
}{}%
\IfNoValueTF{#9}{%
\def\tempc{1}%
}{%
\def\tempc{#9}%
}%

\expandafter\node\expandafter[\tempa](s) at (0,0) {#2};
\expandafter\node\expandafter[\tempb, below = 1*\dl of s](c){#4};
\tio{i1}{s}{east}{1/3};
\tio{i3}{s}{east}{2/3};

\tio{o1}{s}{west}{1/3};
\tio{o3}{s}{west}{2/3};

\draw[->] (o1) -- ($(o1) - (\tempc*\dl, 0)$) |- node[pos = 0.25, left]{#5}(c);
\draw[<-] (i1) -- ($(i1) + (\tempc*\dl, 0)$) |- node[pos = 0.25, right]{#6}(c);

\draw[->] (o3) -- node[at end, above]{#7}($(o3) - (\tempc*\dl + \dl, 0)$);
\draw[<-] (i3) -- node[at end, above]{#8}($(i3) + (\tempc*\dl + \dl, 0)$);
}

\DeclareDocumentCommand \iqccon {o m o m m m m }{%
\def\tempa{#1}%
\def\tempb{#3}%
\IfNoValueTF #1 {%
\def\tempa{sy1}%
}{}%
\IfNoValueTF #3 {%
\def\tempb{sy3}%
}{}%
\expandafter\node\expandafter[\tempa](s) at (0,0) {#2};
\expandafter\node\expandafter[\tempb, above = 1*\dl of s](d){#4};

\tio{i1}{s}{east}{1/2};
\tio{o1}{s}{west}{1/2};

\node[sum, left = 2*\dl of o1](su){$+$};
\draw[->] (o1) -- (su);
\draw[->] (su) |- node[pos=0.25, left]{#5}(d);
\draw[<-] (i1) -- ($(i1) + (2*\dl, 0)$) |- node[pos = 0.25, right]{#6}(d);
\draw[->] ($(su) - (2*\dl, 0)$) -- node[at start]{#7}(su);
}

\DeclareDocumentCommand \iqccona {o m o m m m m m m}{%
\def\tempa{#1}%
\def\tempb{#3}%
\IfNoValueTF #1 {%
\def\tempa{sy1}%
}{}%
\IfNoValueTF #3 {%
\def\tempb{sy3}%
}{}%
\expandafter\node\expandafter[\tempa](s) at (0,0) {#2};
\expandafter\node\expandafter[\tempb, above = 1*\dl of s](d){#4};

\tio{i1}{s}{east}{1/2};
\tio{o1}{s}{west}{1/2};

\tio{i1d}{d}{west}{1/2};
\tio{o1d}{d}{east}{1/2};
\node[sum, left = 2*\dl of o1](su){$+$};
\node[sum, right = 2*\dl of o1d](sub){$+$};

\draw[->] (o1) -- (su);
\draw[->] (su) |- node[pos=0.25, left]{#7}(d);
\draw[->] (d) -- node[above]{#5} (sub);
\draw[->] ($(sub) + (2*\dl, 0)$) -- node[above]{#6} (sub);
\draw[->] (sub) |- node[pos = 0.25,right]{#8}(s);
\draw[->] ($(su) - (2*\dl, 0)$) -- node[at start]{#9}(su);
}

\DeclareDocumentCommand \tcon {o m o m}{%
\def\tempa{#1}%
\def\tempb{#3}%
\IfNoValueTF #1 {%
\def\tempa{sy2}%
}{}%
\IfNoValueTF #3 {%
\def\tempb{sy2}%
}{}%
\expandafter\node\expandafter[\tempa](k) at (0,0) {#2};
\expandafter\node\expandafter[\tempb, right = 2*\dl of k](s){#4};
\tconcontinued
}
\DeclareDocumentCommand \tconcontinued {m m m m m m}{%
\node[sum, left=2*\dl of k](su1){$+$};
\node[sum, right=2*\dl of s](su2){$+$};
\node[sum, right=2*\dl of su2](su3){$+$};
\node[sum, below=2*\dl of su2](su4){$+$};
\draw[->] (su2) -- (su3);
\draw[->] (su2) -- (su4);
\draw[->] (s) -- (su2);
\draw[->] (k) -- node[]{#3} (s);
\draw[->] (su1) -- node[]{#2} (k);
\draw[->] ($(su1) - (3*\dl, 0)$) -- node[at start, above]{#1} (su1);
\coordinate(t1) at ($(su1) - (1.5*\dl, 0)$); 
\draw[->] (t1) -- ($(t1) + (0, 2.5*\dl)$) -| node[pos=0.9]{$-$}(su3);
\draw[->] ($(su2) + (0, 2*\dl)$) -- node[left]{#4} (su2);
\draw[->] ($(su4) + (2*\dl, 0)$) -- node[above]{#6} (su4);
\draw[<-] ($(su3) + (2*\dl, 0)$) -- node[at start, above]{#5}(su3);
\draw[->] (su4) -| node[pos = 0.9]{$-$}(su1);
}

\DeclareDocumentCommand \iosys {o m m m m m}{%
\IfNoValueTF{#1}{%
\expandafter\node\expandafter[#2](sy#6){#3};%
}{
\expandafter\node\expandafter[#2](sy#6)at (#1){#3};%
}%
\foreach \x in {1, ..., #5}{
\tio{i\x#6}{sy#6}{east}{\x/(#5 + 1)};
}
\foreach \x in {1, ..., #4}{
\tio{o\x#6}{sy#6}{west}{\x/(#4 + 1)};
}
}

\DeclareDocumentCommand \lsba {m m m o o}{%
\IfNoValueTF{#5}{%
\def\tempa{1*\dl}%
}{%
\def\tempa{#5}%
}%
\newdimen\mydim
\pgfextractx{\mydim}{\pgfpointdiff{\pgfpointanchor{#1}{east}}{\pgfpointanchor{#2}{east}}}
\ifnum\pdfstrcmp{#4}{end}=0%
	\ifdim\mydim > 0pt%
		\draw[->] (#1) -- ($(#1) - (\tempa, 0)$) |- node[pos = 0.85, above]{#3}(#2);%
	\else%
		\draw[->] (#1) -| ($(#2) - (\tempa, 0)$) -- node[pos = 0.7, above]{#3}(#2);
	\fi
\else%
	\ifnum\pdfstrcmp{#4}{start}=0%
		\ifdim\mydim > 0pt%
			\draw[->] (#1) -- node[pos = 0.3, above]{#3}($(#1) - (\tempa, 0)$) |- (#2);%
		\else%
			\draw[->] (#1) -| node[pos = 0.15, above]{#3}($(#2) - (\tempa, 0)$) -- (#2);
		\fi
	\else%
		\ifdim\mydim > 0pt%
			\draw[->] (#1) -- ($(#1) - (\tempa, 0)$) |- node[pos = 0.25, left]{#3}(#2);%
		\else%
			\draw[->] (#1) -| node[pos = 0.75, left]{#3}($(#2) - (\tempa, 0)$) -- (#2);
		\fi%
	\fi%
\fi%
}

\DeclareDocumentCommand \rsba {m m m o o}{%
\IfNoValueTF{#5}{%
\def\tempa{1*\dl}%
}{%
\def\tempa{#5}%
}%
\newdimen\mydim
\pgfextractx{\mydim}{\pgfpointdiff{\pgfpointanchor{#1}{west}}{\pgfpointanchor{#2}{west}}}
\ifnum\pdfstrcmp{#4}{end}=0%
	\ifdim\mydim < 0pt%
		\draw[->] (#1) -- ($(#1) + (\tempa, 0)$) |- node[pos = 0.85, above]{#3}(#2);%
	\else%
		\draw[->] (#1) -| ($(#2) + (\tempa, 0)$) -- node[pos = 0.7, above]{#3}(#2);%
	\fi%
\else%
	\ifnum\pdfstrcmp{#4}{start}=0%
		\ifdim\mydim < 0pt%
			\draw[->] (#1) -- node[pos = 0.3, above]{#3}($(#1) + (\tempa, 0)$) |- (#2);%
		\else%
			\draw[->] (#1) -| node[pos = 0.11, above]{#3}($(#2) + (\tempa, 0)$) -- (#2);%
		\fi%
	\else%
		\ifdim\mydim < 0pt%
			\draw[->] (#1) -- ($(#1) + (\tempa, 0)$) |- node[pos = 0.25, right]{#3}(#2);%
		\else%
			\draw[->] (#1) -| node[pos = 0.75, right]{#3}($(#2) + (\tempa, 0)$) -- (#2);%
		\fi%
	\fi%
\fi%
}

\DeclareDocumentCommand \bsba {m m m o o}{%
\IfNoValueTF{#5}{%
\def\tempa{1*\dl}%
}{%
\def\tempa{#5}%
}%
\newdimen\mydim
\pgfextracty{\mydim}{\pgfpointdiff{\pgfpointanchor{#1}{west}}{\pgfpointanchor{#2}{west}}}
\ifnum\pdfstrcmp{#4}{end}=0%
		\ifdim\mydim > 0pt%
			\draw[->] (#1) -- ($(#1) - (0,\tempa)$) -|  node[pos= 0.85, left]{#3}(#2);%
		\else%
			\draw[->] (#1) |- ($(#2) -(0, \tempa)$) -- node[pos = 0.7, right]{#3}(#2);%
		\fi%
\else%
	\ifnum\pdfstrcmp{#4}{start}=0%
		\ifdim\mydim > 0pt%
			\draw[->] (#1) -- node[pos= 0.3, right]{#3}($(#1) - (0,\tempa)$) -|  (#2);%
		\else%
			\draw[->] (#1) |- node[pos = 0.15, left]{#3}($(#2) -(0, \tempa)$) -- (#2);%
		\fi%
	\else%
		\ifdim\mydim > 0pt%
			\draw[->] (#1) -- ($(#1) - (0,\tempa)$) -|  node[pos= 0.25, below]{#3}(#2);%
		\else%
			\draw[->] (#1) |- node[pos = 0.75, below]{#3}($(#2) -(0, \tempa)$) -- (#2);%
		\fi%
	\fi%
\fi%
}

\DeclareDocumentCommand \asba {m m m o o}{%
\IfNoValueTF{#5}{%
\def\tempa{1*\dl}%
}{%
\def\tempa{#5}%
}%
\newdimen\mydim
\pgfextracty{\mydim}{\pgfpointdiff{\pgfpointanchor{#1}{west}}{\pgfpointanchor{#2}{west}}}
\ifnum\pdfstrcmp{#4}{end}=0%
		\ifdim\mydim < 0pt%
			\draw[->] (#1) -- ($(#1) + (0,\tempa)$) -| node[pos= 0.85, left]{#3}(#2);%
		\else%
			\draw[->] (#1) |- ($(#2) +(0, \tempa)$) -- node[pos = 0.7, right]{#3}(#2);%
		\fi%
\else%
	\ifnum\pdfstrcmp{#4}{start}=0%
		\ifdim\mydim < 0pt%
			\draw[->] (#1) -- node[pos= 0.3, right]{#3}($(#1) + (0,\tempa)$) -| (#2);%
		\else%
			\draw[->] (#1) |- node[pos = 0.15, left]{#3}($(#2) +(0, \tempa)$) -- (#2);%
		\fi%
	\else%
		\ifdim\mydim < 0pt%
			\draw[->] (#1) -- ($(#1) + (0,\tempa)$) -|  node[pos= 0.25, above]{#3}(#2);%
		\else%
			\draw[->] (#1) |- node[pos = 0.75, above]{#3}($(#2) +(0, \tempa)$) -- (#2);%
		\fi%
	\fi%
\fi%
}

\DeclareDocumentCommand \ioar {m m m m o}{%
\IfNoValueTF{#5}{%
\def\tempa{2*\dl}%
}{%
\def\tempa{#5}%
}%
\ifnum\pdfstrcmp{#3}{north}=0%
	\ifnum\pdfstrcmp{#4}{in}=0%
		\draw[<-] (#1) -- node[right, at end]{#2}($(#1) + (0,\tempa)$);%
	\else%
		\ifnum\pdfstrcmp{#4}{out}=0%
			\draw[->] (#1) -- node[right, at end]{#2}($(#1) + (0,\tempa)$);%
		\fi%
	\fi%
\else%
	\ifnum\pdfstrcmp{#3}{n}=0%
		\ifnum\pdfstrcmp{#4}{in}=0%
			\draw[<-] (#1) -- node[right, at end]{#2}($(#1) + (0,\tempa)$);%
		\else%
			\ifnum\pdfstrcmp{#4}{out}=0%
				\draw[->] (#1) -- node[right, at end]{#2}($(#1) + (0,\tempa)$);%
			\fi%
		\fi%
	\else%
		\ifnum\pdfstrcmp{#3}{south}=0%
			\ifnum\pdfstrcmp{#4}{in}=0%
				\draw[<-] (#1) -- node[right, at end]{#2}($(#1) - (0,\tempa)$);%
			\else%
				\ifnum\pdfstrcmp{#4}{out}=0%
					\draw[->] (#1) -- node[right, at end]{#2}($(#1) - (0,\tempa)$);%
				\fi%
			\fi%
		\else%
			\ifnum\pdfstrcmp{#3}{s}=0%
				\ifnum\pdfstrcmp{#4}{in}=0%
					\draw[<-] (#1) -- node[right, at end]{#2}($(#1) - (0,\tempa)$);%
				\else%
					\ifnum\pdfstrcmp{#4}{out}=0%
						\draw[->] (#1) -- node[right, at end]{#2}($(#1) - (0,\tempa)$);%
					\fi%
				\fi%
			\else%
				\ifnum\pdfstrcmp{#3}{east}=0%
					\ifnum\pdfstrcmp{#4}{in}=0%
						\draw[<-] (#1) -- node[above, at end]{#2}($(#1) + (\tempa,0)$);%
					\else%
						\ifnum\pdfstrcmp{#4}{out}=0%
							\draw[->] (#1) -- node[above, at end]{#2}($(#1) + (\tempa,0)$);%
						\fi%
					\fi%
				\else%
					\ifnum\pdfstrcmp{#3}{e}=0%
						\ifnum\pdfstrcmp{#4}{in}=0%
							\draw[<-] (#1) -- node[above, at end]{#2}($(#1) + (\tempa,0)$);%
						\else%
							\ifnum\pdfstrcmp{#4}{out}=0%
								\draw[->] (#1) -- node[above, at end]{#2}($(#1) + (\tempa,0)$);%
							\fi%
						\fi%
					\else%
						\ifnum\pdfstrcmp{#3}{west}=0%
							\ifnum\pdfstrcmp{#4}{in}=0%
								\draw[<-] (#1) -- node[above, at end]{#2}($(#1) - (\tempa, 0)$);%
							\else%
								\ifnum\pdfstrcmp{#4}{out}=0%
									\draw[->] (#1) -- node[above, at end]{#2}($(#1) - (\tempa, 0)$);%
								\fi%
							\fi%
						\else%
							\ifnum\pdfstrcmp{#3}{w}=0%
								\ifnum\pdfstrcmp{#4}{in}=0%
									\draw[<-] (#1) -- node[above, at end]{#2}($(#1) - (\tempa, 0)$);%
								\else%
									\ifnum\pdfstrcmp{#4}{out}=0%
										\draw[->] (#1) -- node[above, at end]{#2}($(#1) - (\tempa,0)$);%
									\fi%
								\fi%
							\else%
							\fi%
						\fi%
					\fi%
				\fi%
			\fi%
		\fi%
	\fi%
\fi%
}

\newsavebox{\sattikz}
\savebox{\sattikz}{
\begin{tikzpicture}
\draw[->,gray!70!black] (-1, 0) -- (1, 0);
\draw[->,gray!70!black] (0,-0.8) -- (0, 0.8);
\draw[-, rounded corners=false] (-1, -0.5) -- (-0.5, -0.5) -- (0.5, 0.5) -- (1, 0.5);
\end{tikzpicture}}

\DeclareDocumentCommand \satc {o o}{%
\IfNoValueTF{#1}{%
\def\tempa{black}%
}{%
\def\tempa{#1}%
}%
\IfNoValueTF{#2}{%
\resizebox{8ex}{!}{
\begin{tikzpicture}
\draw[->,gray!70!black] (-1, 0) -- (1, 0);
\draw[->,gray!70!black] (0,-0.8) -- (0, 0.8);
\draw[-, rounded corners=false, color=\tempa] (-1, -0.5) -- (-0.5, -0.5) --
(0.5, 0.5) -- (1, 0.5);
\end{tikzpicture}}
}{%
\resizebox{8ex}{!}{
\begin{tikzpicture}
\draw[->,gray!70!black] (-1, 0) -- (1, 0);
\draw[->,gray!70!black] (0,-0.8) -- (0, 0.8);
\draw[-, rounded corners=false, color=\tempa] (-1, -0.5) -- (-0.5, -0.5) --
(0.5, 0.5) -- (1, 0.5);
\draw[-, color=#2] (-0.7, -0.7) -- (0.7, 0.7);
\draw[-, color=#2] (-0.7, 0) -- (0.7, 0);
\end{tikzpicture}}
}%
}

\newsavebox{\deadzonetikz}
\savebox{\deadzonetikz}{
\begin{tikzpicture}
\draw[->, gray!70!black] (-1, 0) -- (1, 0);
\draw[->, gray!70!black] (0,-0.8) -- (0, 0.8);
\draw[-, rounded corners=false] (-1, -0.6) -- (-0.4, 0) -- (0.4, 0) -- (1, 0.6);
\end{tikzpicture}}

\journal{}

\newcommand{\C}{\mathbb C}
\newcommand{\eig}{\text{eig}}
\newcommand{\bul}{\bullet}

\begin{document}

\allowdisplaybreaks

\begin{frontmatter}
\title{On dynamic dissipation inequalities for stability analysis}
\title{Stability analysis by dissipation and dynamic integral quadratic constraints: A complete link}
\title{Stability analysis by dynamic dissipation inequalities: On merging frequency-domain techniques with time-domain conditions}
\author[carsten]{Carsten W. Scherer\fnref{a1}}
\ead{carsten.scherer@mathematik.uni-stuttgart.de}
\author[joost]{Joost Veenman}
\address[carsten]{Department of Mathematics, University of Stuttgart, \\ Pfaffenwaldring 5a,
70569 Stuttgart, Germany}
\address[joost]{Aerospace Division of Sener\\ Severo Ochoa
4 (PTM), 28760, Tres Cantos (Madrid), Spain}
\ead{joostveenman@gmail.com}
\fntext[a1]{
This author would like to thank the German Research Foundation (DFG) for financial support of the project within the Cluster of Excellence in Simulation Technology (EXC 310/2) at the University of Stuttgart.}

\begin{abstract}
In this paper we provide a complete link between dissipation theory and a celebrated result on stability analysis with integral quadratic constraints. This is achieved with a new stability characterization for feedback interconnections based on the notion of  finite-horizon integral quadratic constraints with a terminal cost. As the main benefit, this 
opens up opportunities for guaranteeing constraints on the transient responses of trajectories in feedback loops within absolute stability theory. For parametric robustness, we show
how to generate tight robustly invariant ellipsoids on the basis of a classical frequency-domain stability test, with illustrations by a numerical example.

\end{abstract}
\end{frontmatter}


\section{Introduction}
The framework of integral quadratic constraints (IQCs) was developed in \cite{MegRan97} and builds on the seminal contributions of Yakubovich \cite{Yak67} and Zames \cite{Zam66a,Zam66b}. It provides a technique for analyzing the stability of an interconnection of some linear time-invariant (LTI) system in feedback with another causal system without any particular description, which is also called uncertainty.
The key idea is to capture the properties of the uncertainty through filtered energy relations of the output in response to inputs with finite energy. Mathematically, this is formalized by requiring the $\Le_2$-input-output pairs of the uncertainty to satisfy an integral quadratic constraint, an inequality expressed with a quadratic form on $\Le_2$ in the frequency domain that is defined by a so-called multiplier. In this setting, stability of the interconnection is guaranteed if the LTI system 
satisfies a suitable frequency-domain inequality (FDI) involving the multiplier, which can be computationally verified 
by virtue of the Kalman-Yakubovich-Popov (KYP) lemma. Various papers (cf. \cite{MegRan97,VeeSch16} and references therein) give a detailed exposition of different uncertainties and their corresponding multiplier classes on the basis of which the IQC theorem in \cite{MegRan97} allows to generate a large variety of practical computational robust stability and performance analysis tests.
The stunningly wide impact of this framework also incorporates, e.g., the analysis of adaptive learning \cite{AndYou07} or of optimization algorithms \cite{LesRec16}.

Another central notion in systems theory is dissipativity \cite{Wil72a,Wil72b}, which has been developed by Jan Willems with the explicit goal of arriving at a more fundamental understanding of the stability properties of feedback interconnections \cite[p. 322]{Wil72a}. Roughly speaking, a system with a state-space description is said to be dissipative with respect to some supply rate if there exists a storage function for which a dissipation inequality is valid along all system trajectories; for quadratic supply rates, such dissipation inequalities can be also viewed as integral quadratic constraints.

A huge body of work has been devoted to analyzing the links between both frameworks. In particular if the multipliers (supply rates) are non-dynamic and the two approaches involve so-called hard (finite-horizon) IQCs, the relation between the two worlds is well-established, addressed, e.g., in
\cite{SchWei99,BroLoz07,HadChe08,ArcMei16,Van17};
the classical small-gain, passivity or conic-sector theorems are prominent examples, with generalizations given in \cite{NarTay73,DesVid75,Saf80,Tee96,GeoSmi97}.
However, for the much more powerful dynamic multipliers in \cite{MegRan97}, the connection between the related so-called soft (infinite-horizon) IQCs and dissipation theory has only been demonstrated for specialized 
cases
in  \cite{Bal02,IwaHar05,WilTak07,VeeSch13a,Sei15,CarSei17}. Relations of IQCs to Yakubovich's absolute stability framework and classical multiplier theory are addressed, e.g., in \cite{Bar00,Shi00,Yak00,Yak02,FuDas05,Alt12,CarTur15}.

The purpose of this paper is to present a novel IQC theorem based on the notion of finite-horizon IQCs with a terminal cost. In generalizing \cite{VeeSch13a,Sei15}, a first key contribution is to show that the IQC theorem from \cite{MegRan97} for general multipliers can be subsumed to our framework.
In this way we provide a first ever complete link between the IQC framework and dissipation theory.
Furthermore, we show how a classical frequency domain robust stability test for parametric uncertainties
permits the generation of finite-horizon IQCs with a convexly constrained terminal cost, and illustrate
the ensuing benefit in terms a numerical example. As argued in \cite{Bal02,FetSch17c}, such bridges permit to beneficially merge frequency-domain techniques with time-domain conditions for the construction of local absolute stability criteria, a topic left for future research.


The paper is structured as follows. In Section \ref{S2}, we recall the main IQC theorem and formulate our new dissipation-based stability result. 
Section \ref{S3} develops the relevant technical consequences of the hypotheses in the IQC theorem, which allows to subsume it to our encompassing result in Section \ref{S4}. Finally, Section \ref{S5} illustrates the extra benefit of our
results
over standard IQC theory.

%

{\bf Notation.} $\Le_{2e}^n$ denotes the space of all locally square integrable signals $x:[0,\infty)\to\R^n$, while $\Le_{2}^n:=\{x\in \Le_{2e}^n\st \|x\|:=\sqrt{\int_0^\infty x(t)^Tx(t)\,dt}<\infty\}$. The Fourier transform $x\in\Le_2^n$ is denoted by $\hat x$.
For $T>0$
we work with the truncation operator
$P_T:\Le_{2e}^n\to\Le_{2}^n$, $x\mapsto P_Tx=x_T$ where $x_T=x$ on $[0,T]$ and $x_T=0$ on $(T,\infty)$. The system $S:\Le_{2e}^n\to \Le_{2e}^m$ is casual if $SP_T=P_TSP_T$ for all $T>0$,
and $S$ is bounded (stable) if its induced $\Le_2$-gain is finite.
For a real rational matrix $G$ let $G^*(s):=G(-s)^T$ and
$G=(A,B,C,D)=\begin{scriptsize}\mas{c|c}{A&B\hl C&D}\end{scriptsize}$ means
$G(s)=C(sI-A)^{-1}B+D$. 
Further,  $\rli^{n\times m}$ ($\rhi^{n\times m}$) is the space of real rational matrices without poles on
the extended imaginary axis $\Ci:=i\R\cup\{\infty\}$ (in the closed right-half plane). To save space we use the abbreviations
\equ{\label{KYPg}\arc{0.3}
\kyp{X}{\mi}{\mat{cc}{A&B\\C&D}}:=
\mat{cc}{I&0\\A&B\\C&D}^T
\mat{ccc}{0&X&0\\X&0&0\\0&0&\mi}
\mat{cc}{I&0\\A&B\\C&D}\te{and}
\col(u_1,\ldots,u_n):=\mat{c}{u_1\\\vdots\\u_n}.
}
Finally, ``$\pl$'' stands for objects that can be inferred by symmetry or are irrelevant.

\section{A Novel IQC theorem}\label{S2}

\subsection{Recap of standard IQC theorem}
For setting up the integral quadratic constraint (IQC) framework, we consider the linear finite-dimensional time-invariant system
\eql{sy}{
\arr{rcl}{\dot x&=&\Ag x+\Bg w,\ \ x(0)=0,\\z&=&\Cg x+\Dg w}
}
where $A$ is Hurwitz. This defines a causal linear map $\sy:\Le_{2e}^\nz\to\Le_{2e}^\nw$;
we also denote the related transfer matrix by $\sy(s)=C(sI-A)^{-1}B+D$
which should not cause any confusion. Further, we let $\Delta$ be a system, called uncertainty, which
could be defined through a distributed or nonlinear system but which is not required
to admit any specific (state-space) description. We suppose that
$$\Delta:\Le_{2e}^\nw\to\Le_{2e}^\nz\te{is causal.}$$
In this paper we investigate the feedback interconnection
\eql{fb}{
z=\sy w+d\te{and}w=\Delta(z)
}
of the system $G$ with the uncertainty $\Delta$ that is affected by the external disturbance $d\in\Le_{2e}^\nz$. This loop is said to be
{\em well-posed}  if for any $d\in\Le_{2e}^\nz$ there exists a unique response
$z\in\Le_{2e}^\nz$ that depends causally on $d$; this is equivalent to
$(I-\sy \Delta):\Le_{2e}^\nz\to\Le_{2e}^\nz$ having a causal inverse.
The loop is {\em stable} if there exists some $\ga\geq 0$ such that
$\|z\|\leq \ga\|d\|$ holds for all $d\in L_2^\nz$ and all responses of \r{fb}.
If \r{fb} is well-posed, its stability is equivalent to the inverse $(I-\sy \Delta)^{-1}$ being bounded.


Under the assumptions that $\Delta$ is bounded and the loop \r{fb} is well-posed, the IQC theorem establishes a separation condition on the graph of $\Delta$ and the inverse graph of $\sy$ that guarantees stability of \r{fb}. These conditions
are formulated as frequency domain inequalities (FDIs) for the system $G$ and opposite integral quadratic constraints (IQC) on the uncertainty $\Delta$, both involving some rational transfer matrix $\Pi$, called multiplier, which is assumed to have the properties
$\Pi\in\rli^{(\nz+\nw)\times(\nz+\nw)}\te{and}\Pi=\Pi^*.$
Let us now cite the main theorem of \cite{MegRan97}.
\theorem{\label{Tiqc}Suppose that $\Delta$ is bounded and
that \r{fb} is well-posed for $\tau\Delta$ with any $\tau\in[0,1]$ replacing $\Delta$.
Then
$(I-\sy \Delta)^{-1}$ is bounded if $G$ satisfies the FDI
\eql{fdi}{
\mat{c}{\sy (\io)\\I_\nw}^*\Pi(\io)\mat{c}{\sy (\io)\\I_\nw}\cl 0\te{for all}\io\in\Ci
}
and if for the uncertainty the following IQCs hold:
\eql{iqc}{
\int_{-\infty}^\infty
\mat{c}{\hat z(\io)\\ \tau\widehat{\Delta(z)}(\io) }^*\Pi(\io)
\mat{c}{\hat z(\io)\\ \tau\widehat{\Delta(z)}(\io) }\,d\om\geq 0
\te{for all}z\in\Le_2^\nz,\ \tau\in[0,1].
}
}

\remark{\label{rem1}
We extract from \r{iqc}
for $\tau=0$ that the left-upper $\nz\times\nz$-block of $\Pi$ is positive semi-definite on $\Ci$. If replacing $\Pi$ by $\Pi+\diag(\eps I_\nz,0_\nw)$ for some sufficiently small $\eps>0$, both \r{fdi} and \r{iqc} remain valid. W.l.o.g. we can hence assume in Theorem \ref{Tiqc} that
\eql{lup}{E^T\Pi E\cg 0\te{on}\Ci\te{with}E:=\mat{c}{I_\nz\\0}.}
We say that $\Pi$ is a {\em positive-negative (PN) multiplier} if, next to \r{lup}, its right-lower block is negative semi-definite on the extended imaginary axis:
\eql{rl}{\mat{c}{0\\I_\nw}^T\Pi\mat{c}{0\\I_\nw}\cle 0\te{on}\Ci.}
}

All throughout the paper and without loss of generality $\Pi$ is supposed to be described in terms of a (usually tall) stable outer factor $\Psi$ and a middle matrix $\mi$ as
\eql{fac}{
\Pi=\Psi^*\mi\Psi\te{with}\Psi\in\rhi^{\ny\times(\nz+\nw)}\te{and}\mi=\mi^T\in\R^{n_y\times n_y}.
}
In practice, many multiplier classes do admit the description \r{fac} with some fixed $\Psi$ and a variable matrix $\mi$ (see e.g. \cite{MegRan97,VeeSch16}). If relevant, we work with the state-space description
\eql{psi}{
\arr{rcl}{\dot \xi&=&A_\Psi \xi+B_\Psi u,\ \ \xi(0)=0,\\y&=&C_\Psi \xi+D_\Psi u}
}
of $\Psi$ where $A_\Psi$ is Hurwitz. Viewing \r{psi} as a filter for $u=\col(Gw,w)$ or $u=\col(z,\Delta(z))$ allows to translate the FDI and IQC in Theorem \ref{Tiqc} into time-domain
dissipation inequalities as seen next.

%

\subsection{Main result}\label{Stro}

On the basis of \r{sy} and \r{psi} let us introduce the realizations
\eql{psisy}{
F:=\mat{c}{G\\I}=\mas{c|c}{\Ag&\Bg\hl\Cg&\Dg\\0&I}=:\mas{c|c}{\Ag&\Bg\hl\Cgr&\Dgr},\ \
\Psi\mat{c}{\sy\\I}=
\mas{cc|c}{A_\Psi&B_\Psi\Cgr&B_\Psi \Dgr\\0&\Ag&\Bg\hl C_\Psi&D_\Psi\Cgr& D_\Psi\Dgr}=:
\mas{c|c}{\Af&\Bf\hl \Cf&\De}
}
for the transfer matrix of the system's inverse graph and the filtered version thereof.
Since $A$ is Hurwitz and by the KYP-Lemma, \r{fdi} holds iff there exists some $\Xf=\Xf^T$ with
\eql{kyp}{
\kyp{\Xf}{\mi}{\mat{cc}{\Af&\Bf\\\Cf&\Df}}\cl 0;
}
in the sequel we say that $\Xf$ or \r{kyp} certify the FDI \r{fdi}, or that $\Xf$ is a certificate thereof; moreover, whenever relevant we assume $\Xf$ to be partitioned as $\Af$ in \r{psisy}.
%
%

Now suppose that \r{kyp} is valid. By Finsler's lemma, we can choose some $\ga>0$ with
\eql{lmip}{
\kyp{\Xf
}{
\mat{ccccc}{
\mi&0&0\\
0&\frac{1}{\ga}I&0\\
0&0             &-\ga I}}{
\mat{cc|cccc}{
A_\Psi&B_\Psi\Cgr&
B_\Psi\Dgr&
B_\Psi \Ddgr\\0&A&B&0\hl
C_\Psi&D_\Psi\Cgr& D_\Psi\Dgr&D_\Psi \Ddgr\\
0&\Cg&\Dg&I_{\nz}\\
0&0&0&I_{\nz}}}
\cl 0.
}
This step leads to a crucial dissipation inequality as follows. If $z=Gw+d$ is the response to any $w\in\Le_{2e}^\nw$ and $d\in\Le_{2e}^\nz$ and if we let $u=\col(z,w)$ drive the filter \r{psi}, we have
$$
\arr{rcl}{
\mat{c}{\dot \xi\\\dot x}&=&
\mat{cc}{
A_\Psi&B_\Psi\Cgr\\
0     &A         }
\mat{c}{\xi\\x}+
\mat{cc}{
B_\Psi\Dgr& B_\Psi\Ddgr\\
B&0}\mat{c}{w\\d},\ \ \mat{c}{\xi(0)\\x(0)}=0\\
\mat{c}{y\\z\\d}&=&
\mat{cc}{
C_\Psi&D_\Psi\Cgr\\
0&\Cg\\
0&0}\mat{c}{\xi\\x}+
\mat{cccccc}{
D_\Psi\Dgr&D_\Psi \Ddgr\\
\Dg&I_\nz\\
0&I_\nz}\mat{c}{w\\d}.
}$$
With the combined state trajectory $\eta=\col(\xi,x)$, we right- and left-multiply
\r{lmip} by $\col(\eta,w,d)$ and its transpose to obtain
$$
\frac{d}{dt}\eta(t)^T\Xf\eta(t)+y(t)^TMy(t)+\frac{1}{\ga}\|z(t)\|^2-\ga\|d(t)\|^2\leq 0\te{for}t\geq 0.
$$
After integration we arrive at the dissipation inequality
\eql{dis}{
\eta(T)^T\!\!\Xf\eta(T)+\!\! \int_0^T\!\! y(t)^T \mi y(t) \, dt+\!
\int_0^T\! \frac1\ga \|z(t)\|^2 -\ga \|d(t)\|^2 \, dt
\leq 0\te{for}T>0.
}

On the other hand, let us consider the IQC \r{iqc} for $\tau=1$ and if $\Delta$ is bounded:
\eql{iqc1}{
\int_{-\infty}^\infty
\mat{c}{\hat z(\io)\\ \widehat{\Delta(z)}(\io) }^*\Psi(\io)^*\mi\Psi(\io)
\mat{c}{\hat z(\io)\\ \widehat{\Delta(z)}(\io) }\,d\om\geq 0
\te{for all}z\in\Le_2^\nz.
}
For $z\in\Le_{2}^\nz$ and $u=\col(z,\Delta(z))$ driving \r{psi}, Parseval's theorem shows
that \r{iqc1} implies the validity of the so-called soft (infinite-horizon) IQC
\eql{iqctd}{\int_0^\infty y(t)^T\mi y(t)\,dt\geq 0.}


If $z$ in response to $d\in\Le_2^\nz$ in \r{fb} has finite energy, one can easily infer stability of \r{fb} as follows: Due to $\int_0^T y(t)^T \mi y(t) \, dt\to\int_0^\infty y(t)^T \mi y(t) \, dt$
and $\eta(T)\to 0$ for $T\to\infty$, we can just combine \r{dis} with \r{iqctd} to get $\|z\|\leq\ga\|d\|$. The key difficulty is to {\em prove} that $d\in\Le_2^\nz$ in \r{fb} implies $z\in\Le_{2}^\nz$. This is simple if $\Xf$ in \r{dis} is positive definite and if $z\in\Le_{2e}^\nz$ with $u=\col(z,\Delta(z))$ driving \r{psi} leads to validity of the so-called hard (finite-horizon) IQC
\eql{iqch}{\int_0^T y(t)^T\mi y(t)\,dt\geq 0\te{for}T>0.}
Then \r{dis} guarantees $\|z_T\|\leq\ga\|d_T\|$ for all $T>0$, and $d\in\Le_2^\nz$ implies $\|z\|\leq\ga\|d\|$ by taking the limit $T\to\infty$.
In general, however, neither does \r{kyp} have a positive definite solution, nor can one replace \r{iqctd} with the hard IQC \r{iqch} \cite{MegRan97,VeeSch13a,Sei15}.

Instead of hard and soft IQCs, we propose the following seemingly new notion.

\definition{The uncertainty $\Delta$ satisfies a finite-horizon IQC with terminal cost matrix $Z=Z^T$
(and with respect to the factorization $\Pi=\Psi^*\mi\Psi$) if
\eql{disu}{
\int_0^T y(t)^T\mi y(t)\,dt+\xi(T)^T Z\xi(T)\geq 0\te{for}T>0
}
holds for the trajectories of the filter \r{psi} driven by $u=\col(z,\Delta(z))$ with any
$z\in\Le_{2e}^\nz$.
}

Based on the above line of reasoning and the correct positivity hypothesis on certificates, the following main result of this paper has a simple proof.


\theorem{\label{Tmain}
Let $\Delta$ satisfy a finite-horizon IQC with terminal cost matrix $Z=Z^T$ and suppose
that \r{kyp} has a solution $\Xf=\Xf^T$ which is coupled with $Z$ as
\eql{pos}{
\mat{cc}{\Xf_{11}-Z&\Xf_{12}\\\Xf_{21}&\Xf_{22}}\cg 0.
}
Then there exists some $\ga>0$ such
\eql{disfb}{
\mat{c}{\xi(T)\\x(T)}^T\!\!\mat{cc}{\Xf_{11}-Z&\Xf_{12}\\\Xf_{21}&\Xf_{22}} \mat{c}{\xi(T)\\x(T)}+\!
\int_0^T\! \frac1\ga \|z(t)\|^2 -\ga \|d(t)\|^2 \, dt
\leq 0\te{for}T>0}
holds along the trajectory of the filter \r{psi} driven by any response $u=\col(z,w)$
of the feedback interconnection \r{fb} to any disturbance $d\in\Le_{2e}^\nz$. Moreover,
$d\in\Le_2^\nz$ implies $z\in\Le_2^\nz$ as well as $\|z\|\leq\ga\|d\|$ for all responses of \r{fb}.
}


{\bf Proof.} Let $\Xf$ satisfy \r{kyp} and choose $\ga>0$ such that \r{dis} is valid along the filtered trajectories of \r{fb}.  By assumption, we infer from $w=\Delta(z)$ that
\r{disu} is valid. If just subtracting  \r{disu} from \r{dis} we obtain \r{disfb}.
Now let $d\in\Le_2^\nz$. Since $\Xf-\diag(Z,0)>0$, we infer from \r{disfb} that $\frac{1}{\ga}\|z_T\|^2\leq\ga\|d_T\|^2\leq\ga\|d\|^2$ for all $T>0$, which implies $z\in\Le_{2}^\nz$ and
$\|z\|\leq\ga\|d\|$ by taking the limit $T\to\infty$. 
\epro

We emphasize that Theorem \ref{Tmain} neither requires $G$ or $\Delta$ to be stable nor \r{fb} to be well-posed. Still, it even provides sharper conclusions about the responses of \r{fb} than mere stability, since \r{disfb}
provides information about hard ellipsoidal time-domain constraints on the transient behavior of the
filtered system's state-trajectory in the feedback loop!

A large variety of stability results can be subsumed to Theorem \ref{Tmain}. As one of the core technical contributions of this paper, we reveal that this holds true for Theorem \ref{Tiqc} and general multipliers. For this purpose, we show that \r{fdi} implies the existence of a solution of \r{kyp} with \r{pos} for some suitable matrix $Z$. As a next step, we prove that $\Delta$ satisfies a finite-horizon IQC with a terminal cost for the very same matrix $Z$. Taken together, we conclude that the hypotheses of Theorem \ref{Tiqc} imply that those of Theorem \ref{Tmain} are valid. This not only unveils an unprecedented dissipation proof of Theorem \ref{Tiqc}, but it also allows to draw all the
conclusions for the transient behavior of trajectories in Theorem \ref{Tmain} under the assumptions of the general IQC Theorem \ref{Tiqc}.





\section{A dissipation proof of the IQC theorem for positive negative multipliers}\label{S3}

\subsection{On canonical factorizations}\label{Scf}

It has already been observed in \cite{Sei15,VeeSch13a} that an important role in a dissipation proof of the IQC theorem is played by replacing
\r{fac} with a so-called canonical (Wiener-Hopf) or $J$-spectral factorization
\eql{js}{
\Psi^*\mi\Psi=\t\Psi^*\t\mi\t\Psi\te{with real\ \ $\t \mi=\t \mi^T$\ \ and}\t\Psi,\ \t\Psi^{-1}\in\rhi^{(\nz+\nw)\times(\nz+\nw)}.}
Since $\t\Psi(\infty)$ is non-singular and in view of \r{fdi} and Remark \ref{rem1}, $\t M$ has $\nz$ positive and $\nw$ negative eigenvalues; it can even be chosen as $\diag(I_\nz,-I_\nw)$, but this is not essential in the current paper.

It is well-known \cite{Mei95} that a factorization \r{js} exists if $D_\Psi^T\mi D_\Psi$ is invertible
and if the following algebraic Riccati equation (ARE) has a (unique) stabilizing solution $\t Z=\t Z^T$:
\eql{are}{
A_\Psi^T\t Z+\t ZA_\Psi+C_\Psi^T\mi C_\Psi-(\t ZB_\Psi+C_\Psi^T\mi D_\Psi)(D_\Psi^T\mi D_\Psi)^{-1}(B_\Psi^T\t Z+D_\Psi^T\mi C_\Psi)=0.}
With $\t \mi :=D_\Psi^T\mi D_\Psi$, $\t C_\Psi:=\t \mi^{-1}(B_\Psi^T\t Z+D_\Psi^T\mi C_\Psi)$ and $\t D_\Psi:=I$,
the ARE \r{are} implies
\eql{are1}{
\mat{cc}{I&0\\A_\Psi&B_\Psi\\C_\Psi&D_\Psi}^T\mat{ccc}{0&\t Z&0\\\t Z&0&0\\0&0&\mi }\mat{cc}{I&0\\A_\Psi&B_\Psi\\C_\Psi&D_\Psi}=
\mat{cc}{\t C_\Psi&\t D_\Psi}^T\t \mi \mat{cc}{\t C_\Psi&\t D_\Psi},
}
and the fact that $\t Z$ is stabilizing translates into
\eql{are1s}{\eig(A_\Psi-B_\Psi\t D_\Psi^{-1}\t C_\Psi)\subset\Cl.}
Conversely, if $\det(\t D_\Psi)\neq 0$ then \r{are1}-\r{are1s} imply that $\t Z$ is the stabilizing solution of \r{are}.

Right-multiplying \r{are1} with $\col(\io I-A_\Psi)^{-1}B_\Psi,I)$ and left-multiplying the conjugate transpose indeed shows, after a routine computation, that
\r{are1}-\r{are1s} leads to \r{js} for $\t \Psi=(A_\Psi,B_\Psi,\t C_\Psi,\t D_\Psi)$.
We express this fact by saying that $\t Z$ is a certificate, or certifies \r{js}.
Next we show that $\t Z$ is a candidate for embedding Theorem \ref{Tiqc} into our main result.

\subsection{Consequences of the KYP inequality for the system $G$}\label{section3p2}

Recall that \r{fdi} implies the existence of a certificate
$\Xf=\Xf^T$
for
\r{kyp}. Moreover, the additional property \r{lup} guarantees  the existence of a canonical factorization of $\Pi=\Psi^*\mi\Psi$ as certified by some $\t Z$ \cite{VeeSch13a}. This is stated in the first part of the following result, while the second part assures an instrumental property of positivity.

\lemma{\label{Lfdi2}
Suppose that \r{fdi} holds and that the multiplier satisfies \r{lup}. Then \r{are} has a stabilizing solution $\t Z$. Moreover, $\Xf-\diag(\t Z,0)\cg 0$ holds for all $\Xf=\Xf^T$ with \r{kyp}.
}
%

The proof is found in \ref{apndx2}. The line of reasoning is as follows.
As the first step, we diagonally combine the two FDIs \r{lup} and \r{fdi} with $F$ in \r{psisy}
to get
\eql{fdic}{
\mat{cc}{\t\Psi E&0\\ 0&\t\Psi F}^*
\mat{cc}{-\t\mi&0\\0&\t\mi}
\mat{cc}{\t\Psi E&0\\ 0&\t\Psi F}\cl 0\te{on}\Ci.}
Now observe that both $\mat{cc}{E&F}$ and $\mat{cc}{E&F}^{-1}$ are stable, which implies the same for
$H:=\t \Psi\mat{cc}{-E&F}\mat{cc}{E&F}^{-1}\t \Psi^{-1}$.
With $\small T:=\mat{cc}{-I_{\nz+\nw}&I\\I&I_{\nz+\nw}}$,
a simple computation shows
\eql{tr1}{
\mat{cc}{-\t\mi&0\\0&\t\mi}=T^T\mat{cc}{0&\frac{1}{2}\t\mi\\\frac{1}{2}\t\mi&0}T\te{and}
T\mat{cc}{\t\Psi E&0\\ 0&\t\Psi F}=
\mat{cc}{H\\ I}\t \Psi\mat{cc}{E&F}.
}
With \r{tr1} and since $\t \Psi\mat{cc}{E&F}$ is invertible on $\Ci$, we conclude that \r{fdic} is equivalent to
\eql{fdit}{
\mat{cc}{H\\I}^*
\mat{cc}{0&\frac{1}{2}\t\mi\\\frac{1}{2}\t\mi&0}
\mat{cc}{H\\I}\cl 0\te{on}\Ci.
}
The idea of the proof is to follow these steps for the corresponding KYP inequalities.
With certificates for \r{lup} and \r{fdi}, we construct a certificate for \r{fdit}, which
must be positive definite due to stability of $H$ and ``positive real'' structure of the FDI.
This allows to extract the claimed positivity property.

In \cite{Sei15}, a proof of Lemma \ref{Lfdi2} is given by game-theoretic arguments
(involving both the system FDI and the uncertainty IQC) and by making essential use of the extra constraint \r{rl}. For general multipliers, Lemma \ref{Lfdi2} and the technique of proof are new.

\subsection{Consequences of the IQC for the uncertainty $\Delta$}

If $\Pi=\Psi^*\mi\Psi$ admits a canonical factorization, we now establish
that \r{iqc1} implies the validity of a finite-horizon IQC. This is first formulated
for multipliers satisfying \r{rl}.

\lemma{\label{Liqc}
Let $\Delta$ be stable and satisfy \r{iqc1}.
Moreover, suppose that $\Pi=\Psi^*\mi\Psi$ with \r{rl} admits a canonical factorization as certified by $\t Z$. Then $\Delta$ satisfies a finite-horizon IQC with terminal cost matrix $\t Z$.
}

The proof is found in \ref{apndx3}.
This generalizes \cite{Sei15}, where it is shown that $\Delta$ satisfies a finite-horizon IQC with terminal cost matrix $0$ w.r.t.\  the canonical factorization $\Pi=\t\Psi^*\t\mi\t\Psi$.
Our more general version applies to any multiplier factorization $\Pi=\Psi^*\mi\Psi$
(as a long as it admits a canonical one) and admits a useful extension as seen in Section~\ref{S4}.

\subsection{The IQC theorem for PN multipliers}\label{Sppn}


The following result is a consequence of combining Lemmas \ref{Lfdi2}--\ref{Liqc} and Theorem \ref{Tmain}.

\theorem{\label{Lcom}
Let $\Pi=\Psi^*\mi\Psi$ have the properties \r{lup} and \r{rl}. If
$G$ satisfies the FDI \r{fdi}, $\Delta$ is stable and fulfills \r{iqc1}, then
the hypotheses of Theorem \ref{Tmain} are fulfilled for $Z:=\t Z$.}

If the assumptions in Theorem \ref{Lcom} are valid, one can guarantee
$$
\ga\mat{c}{\xi(T)\\x(T)}^T\!\!
\mat{cc}{\Xf_{11}-\t Z&\Xf_{12}\\\Xf_{21}&\Xf_{22}}\mat{c}{\xi(T)\\x(T)}+\!
\int_0^T\! \|z(t)\|^2\,dt\leq \ga^2\int_0^T\! \|d(t)\|^2\,dt\te{for}T>0
$$
along any loop trajectory that drives the filter \r{psi}.
If $d\in\Le_2^{\nz}$ and in view of \r{pos} for $Z:=\t Z$, this leads to
an ellipsoidal bound on $\col(\xi,x)$ and an energy bound on $z$ in terms of the energy of the disturbance input $d$, even if \r{fb} is not well-posed.
For a rather elaborate discussion on the consequences of such bounds in the IQC setting, we refer to \cite{FetSch17c}.

If starting with a different initial multiplier description $\Pi=\Psi_0^*\mi_0\Psi_0$, the same conclusions can be drawn with a certificate $\t Z_0$ for a corresponding canonical factorization
$\Pi=\t \Psi_0^*\t \mi_0\t \Psi_0$. In case that $\Pi=\Psi_0^*\mi_0\Psi_0$ itself already is a canonical factorization, we get $\t Z_0=0$ and Theorem \ref{Lcom} recovers the main stability result in \cite{Sei15}. Our approach has the benefit that it can be applied to any initial factorization $\Pi=\Psi^*\mi\Psi$ and clearly exhibits the
conceptual role of the certificate $\t Z$ in the conclusions.

\section{A dissipation proof of the general IQC theorem}\label{S4}

To date, no dissipation proofs exist for multipliers which do not satisfy \r{rl}, which excludes several important classes of practical relevance \cite{MegRan97,FetSch17a}. In moving towards overcoming this deficiency,
we show that Lemma \ref{Liqc} persists to hold if replacing \r{rl} with
\eql{grl}{
\mat{c}{H\\I}^*\Pi\mat{c}{H\\I}\cle 0
\te{on}\Ci}
with a stable transfer matrix $H$ for which $z=Hw+d$, $w=\Delta(z)$ is well-posed and stable.

\lemma{\label{Liqc2}
Suppose that $\Delta$ is bounded and satisfies \r{iqc1}.
Moreover, let $\Pi=\Psi^*\mi\Psi$ admit a canonical factorization \r{js} as certified by $\t Z$
and let \r{grl} hold for some $H\in\rhi^{\nz\times\nw}$ for which $(I-H\Delta)^{-1}$ is causal and bounded. Then $\Delta$ satisfies a finite-horizon IQC with terminal cost matrix $\t Z$.
}
%
The proof is given in \ref{apndx4}. This result is new and encompasses Lemma \ref{Liqc} with the choice $H=0$. It is the key technical step in proving the following embedding result.



\theorem{\label{Tnew}Under the hypotheses of Theorem \ref{Tiqc}, the multiplier
$\Pi=\Psi^*\mi\Psi$ with \r{lup} admits a canonical factorization with certificate $\t Z=\t Z^T$,
and the hypotheses of Theorem \ref{Tmain} are fulfilled for $Z:=\t Z$.}
{\bf Proof.} In view of \r{lup}, Lemma \ref{Lfdi2} guarantees the existence of $\t Z=\t Z$ with the property
$\Xf-\diag(\t Z,0)\cg 0$ for any certificate $\Xf=\Xf^T$ of \r{fdi}.
It remains to show that $\Delta$ satisfies a finite-horizon IQC with terminal cost matrix $\t Z$.

The idea is to exploit \r{fdi} which allows us to choose some $\rho\in (0,1)$ (close to $1$) with
\eql{rlt}{\mat{c}{\rho\sy\\I}^*\Pi\mat{c}{\rho\sy\\I}\cl 0\te{on}\Ci.}
Since $\rho^k\to 0$ for $k\to\infty$ and by the small-gain theorem, we can fix some non-negative  integer $k_0$ such that $(I-\rho^{k_0+1}\sy\Delta)^{-1}$ is bounded. Since the latter equals
$(I-[\rho \sy][\rho^{k_0}\Delta])^{-1}$ and with \r{rlt}, the hypotheses of Lemma \ref{Liqc2} hold for $H:=\rho\sy$ and $\rho^{k_0}\Delta$; therefore, $\rho^{k_0}\Delta$ satisfies a finite-horizon IQC with terminal cost matrix $\t Z$.

Now suppose this to be true for $\rho^{k}\Delta$ for any positive integer $k$. By Theorem \ref{Tmain}
applied to the system $\sy$ and the uncertainty $\rho^{k}\Delta$, we infer that
$(I-\rho^{k}\sy\Delta)^{-1}$ is bounded. As just argued, this in turn shows that
$\rho^{k-1}\Delta$ satisfies a finite-horizon IQC with terminal cost matrix $\t Z$.
By induction, this property stays true for all $k=k_0,k=k_0-1,\ldots,k=1$.
\epro

In summary, the hypotheses in Theorem \ref{Tiqc} allow to draw exactly the same conclusions
as for Theorem \ref{Lcom} in Section \ref{Sppn}. In particular the possibility to infer
input-to-state properties from the assumptions in Theorem \ref{Tiqc} and for general multiplies is new.
We emphasize that the choice $\t Z$ is just one 
of many possibilities
to achieve the embedding of Theorem~\ref{Tiqc} into Theorem \ref{Tmain}. In particular in view of computational aspects \cite{FetSch17c},
it is promising to explore in how far one can directly apply Theorem \ref{Tmain} with
matrices $Z=Z^T$ that are not constrained by a (non-convex) ARE but that vary in some well-specified convex set. As another contribution of this paper, one such instance is revealed next.



\section{An application}\label{S5}
The subsequent example serves to illustrate the benefits of Theorem~\ref{Tmain} over Theorem~\ref{Tiqc}. We assume $\nw=\nz$ and consider the class $\bm{\Delta}$
of multiplication operators $\Delta_\delta:\Le_{2e}^\nz\to\Le_{2e}^\nz$, $\Delta_\delta(w)(t):=\delta w(t)$ for $t\geq 0$, with an arbitrary $\delta\in[\al,\be]$ and     for fixed $\al,\be\in\R$ with $0\in(\al,\be)$.
The loop \r{fb} is well-posed for all uncertainties in $\bm{\Delta}$ iff $\det(I-D\delta)\neq 0$ for all $\delta\in[\al,\be]$, which is assumed from now on.

To continue, recall that a transfer matrix $H$ is called generalized positive real (GPR) if $H\in\rli^{\nz\times \nz}$ and if it satisfies $H^*+H\cg 0$ on $\Ci$.
Then the following frequency domain stability characterization essentially goes back to \cite{BroWil65,NarTay73}.
In the terminology of \cite{PacDoy93,ZhoDoy96}, this means that
the verification of robust stability of a feedback interconnection involving one real repeated parametric uncertainty block with dynamic $D/G$ scalings is exact.

\theorem{\label{TBW}The loop \r{fb} is stable for all $\Delta\in\bm{\Delta}$ if and only if there exists some $H$ which is GPR such that $(G-\be^{-1}I)^*H(\al G-I)$ is GPR.}

In order to obtain a computational test, we fix two stable transfer matrices $\psi_1\in\rhi^{m_1\times\nz}$, $\psi_2\in\rhi^{m_2\times\nz}$ and search $H$ in Theorem \ref{TBW} among all transfer matrices $\psi_1^*P\psi_2$ with a free $P\in\R^{m_1\times m_2}$ to render
$\psi_1^*P\psi_2$ and $-(G-\be^{-1}I)^*\psi_1^*P\psi_2(-\al G+I)$ GPR.
With
$$
\Psi:=\mat{cc}{\psi_1&0\\0&\psi_2},\ \ J:=\mat{c}{I\\I},\ \ T:=\mat{cc}{I&-\be^{-1}I\\-\al I&I}
\te{and}
\c M:=\left\{\mat{cc}{0&P\\P^T&0}\st P\in\R^{m_1\times m_2}\right\},
$$
this precisely amounts to testing whether there exists some $\mi\in\c M$ such that
\addtocounter{equation}{1}
\edef\pora{(\arabic{equation}a)}
\edef\porb{(\arabic{equation}b)}
\eql{por}{
\tag{\arabic{equation}a,b}[\Psi J]^*M[\Psi J]\cg 0\te{and}\mat{c}{G\\I}^*[\Psi T]^*M[\Psi T]\mat{c}{G\\I}\cl 0\te{on}\Ci.}


\remark{\label{rem2} Let us take, e.g.,
$\psi_{1}(s)=\psi_2(s)=\col\left(1,\frac{1}{s+1},\ldots,\frac{1}{(s+1)^\nu}\right)$ for $\nu=0,1,2,\ldots$
in \r{por}. If these FDIs are valid for some $\nu\in\N_0$ and $M\in\c{M}$, then robust stability is guaranteed by Theorem \ref{TBW}.
Now recall that any $H$ which is GPR can be uniformly approximated on $\Ci$ by $\psi_1^*P\psi_2$ with a suitable real matrix $P$ (for sufficiently large $\nu$) \cite[Lemma 6]{Sch13}. By Theorem
\ref{TBW}, we can hence conclude that robust stability guarantees the existence of
$\nu\in\N_0$ and $M\in\c{M}$ for which \r{por} is valid; in this sense the proposed robust stability test is asymptotically (i.e, for $\nu\to\infty$) exact.
}

For the multiplier $\Pi=[\Psi T]^*M[\Psi T]$ and any $\delta\in[\al,\be]$, we note that
\eql{fdiu}{
\mat{c}{I\\\delta I}^T\Pi\mat{c}{I\\\delta I}=
[\Psi J]^*M[\Psi J](1-\be^{-1}\delta)(\delta-\al)\cge 0\te{on}\Ci,
}
which implies that $\Delta_\delta$ satisfies \r{iqc1}. For $\delta=0\in(\al,\be)$, \r{fdiu} is strict and $\Pi$ thus satisfies \r{rl}. Since $\Pi's$ right-lower block equals $-\be^{-1}[\Psi J]^*M[\Psi J]$, it actually is a PN-multiplier.

With minimal state-space realizations $\psi_j=(A_j,B_j,C_j,D_j)$ for $j=1,2$, we get with
$(A_\Psi,B_\Psi,C_\Psi,D_\Psi):=(\diag(A_1,A_2),\diag(B_1,B_2),\diag(C_1,C_2),\diag(D_1,D_2))$
a realization of $\Psi$, and the FDIs \r{por} are certified by
\addtocounter{equation}{1}
\edef\lmiRa{(\arabic{equation}a)}
\edef\lmiRb{(\arabic{equation}b)}
\eql{lmiR}{\kyp{R}{M}{\mat{c|c}{A_\Psi&B_\Psi J\hl C_\Psi&D_\Psi J}}\cg 0
\te{and}
\kyp{\Xf}{\mi}{
\mat{cc|c}{A_\Psi&B_\Psi T\Cgr&B_\Psi T \Dgr\\0&\Ag&\Bg\hl C_\Psi&D_\Psi T\Cgr& D_\Psi T\Dgr}}\cl 0.
\tag{\arabic{equation}a,b}
}
Note that $\Xf$ now carries a $3\times 3$ partition.

\lemma{\label{LK}Suppose \r{lmiR} hold for $M\in\c{M}$. Then $\Pi:=[\Psi T]^*M[\Psi T]$ has a canonical factorization that is certified by some $\t Z$ which admits the structure
\eql{cer}{\small\mat{cc}{0&K\\K^T&0}.}
Moreover, any $\Delta\in\bm{\Delta}$ satisfies a finite
horizon IQC with terminal cost matrix $\t Z$ and
\addtocounter{equation}{1}
\eql{cou}{\tag{\arabic{equation}a,b}
\mat{ccc|ccc|cc}{
R_{11} &R_{12}-K     \\
R_{21}-K^T &R_{22}   }\cl 0
\te{as well as}
\mat{cccccccc}{
\Xf_{11} &\Xf_{12}-K &\Xf_{13}    \\
\Xf_{21}-K^T &\Xf_{22} &\Xf_{23}  \\
\Xf_{31} &\Xf_{32} &\Xf_{33}  }\cg 0.
}
}
\edef\coua{(\arabic{equation}a)}
\edef\coub{(\arabic{equation}b)}
In view of Lemmas \ref{Lfdi2}--\ref{Liqc}, we only need to argue why $\t Z$ has the structure
\r{cer} and satisfies $R-\t Z\cl 0$ as in \coub. This is done in \ref{apndx5}.

If there exist $M\in\c{M}$, $\Xf=\Xf^T$, $R=R^T$ with \r{lmiR} and \coub, Theorem~\ref{Tmain}
guarantees robust stability of \r{fb} for all $\Delta\in\bm{\Delta}$, with all consequences on ellipsoidal invariance based on the matrix \coub\ as discussed in Section~\ref{Sppn}. In a computational stability test, it is desirable to view all $\mi$, $\Xf$, $R$ and, in particular, $K$ as decision variables in a convex program. However, this is prevented by the fact that $K$ needs to satisfy an indefinite ARE
(see the proof of Lemma \ref{LK}). This leads us to the last contribution of this paper. We can shown
that it is possible to replace the non-convex ARE constraint on $K$ by the convex constraint \coua\
and still obtain guarantees for robust stability.

\theorem{\label{Tpar}Let there exist $M\in\c M$, $R=R^T$ and $K$ with \lmiRa\ and \coua.
Then any $\Delta\in\bm{\Delta}$ satisfies a finite horizon IQC with terminal cost matrix \r{cer}.}

{\bf Proof.} Let $y$ denote the response of \r{psi} for $u=\col(u_1,u_2)\in\Le_{2e}^{\nz+\nz}$. If $u_1=u_2$ then \lmiRa\ implies $\int_0^T y(t)^TMy(t)\,dt+ \xi(T)^TR\xi(T)\geq 0$ for $T>0$;
combined with \coua, we get
\eql{fhi}{
\int_0^T \mat{c}{y_1(t)\\y_2(t)}^TM\mat{c}{y_1(t)\\y_2(t)}\,dt+ \mat{c}{\xi_1(T)\\\xi_2(T)}^T\mat{cc}{0&K\\K^T&0}\mat{c}{\xi_1(T)\\\xi_2(T)}\cge 0\te{for}T>0.
}

Now let $\delta\in[\al,\be]$. For any $z\in\Le_{2e}^\nz$ set $w=\Delta_\delta(z)=\delta z$ and consider $y=[\Psi T]\col(z,w)$ which equals the response of the filter \r{psi} for
$$
\mat{c}{u_1\\u_2}=T\mat{c}{z\\\delta z}=\mat{c}{(\be-\delta)z\\(\delta-\al)z}.
$$
We claim that \r{fhi} persists to hold for this trajectory. If $\delta=\al$ (or $\delta=\be$), this is trivial since then $u_2=0$ and thus $\xi_2=0$, $y_2=0$ (or $u_1=0$ and thus $\xi_1=0$, $y_1=0$).
If $\delta\in(\al,\be)$, define $\t\delta:=(\delta-\al)/(\be-\delta)\in(0,\infty)$ to infer $u_2=\t\delta u_1$; by linearity, $\col(\t\delta\xi_1,\xi_2)$ and
$\col(\t\delta y_1,y_2)$ are the state- and output responses of the filter to $\col(\t\delta u_1,\t\delta
u_1)$; hence \r{fhi} shows
\eqn{
\int_0^T \mat{c}{\t \delta y_1(t)\\y_2(t)}^T\mat{cc}{0&P\\P^T&0}\mat{c}{\t \delta y_1(t)\\y_2(t)}\,dt+ \mat{c}{\t \delta\xi_1(T)\\\xi_2(T)}^T\mat{cc}{0&K\\K^T&0}\mat{c}{\t \delta\xi_1(T)\\\xi_2(T)}\cge 0
\te{for}T>0.
}
The particular structure allows us to divide by $\t\delta>0$, which indeed leads to \r{fhi}.
\epro

In summary, if there exist $M\in\c M$, $R=R^T$, $\Xf=\Xf^T$ and $K$ with \r{lmiR} and \r{cou},
then the hypothesis of Theorem \ref{Tmain} are satisfied for all $\Delta\in\bm{\Delta}$ and \r{fb} is robustly stable. Following \cite{Bal02,FetSch17c}, it is now routine to proceed as in the following numerical example.

\begin{figure}[bh]\center
\includegraphics[width=13cm]{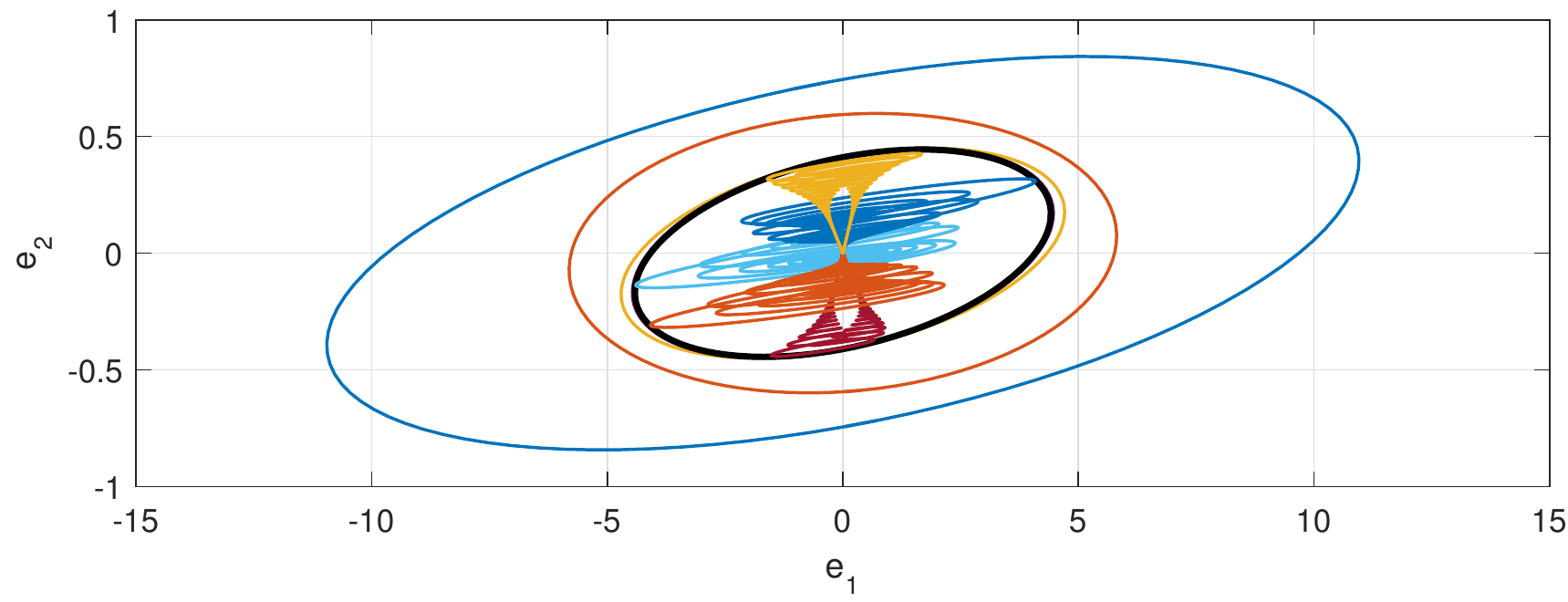}
\caption{Ellipsoids computed for $\nu=0$ (blue), $\nu=1$ (red), $\nu=2$ (yellow), $\nu=3$ (black), and several worst-case system trajectories with disturbances of energy one.}\label{fig1}
\end{figure}

{\bf Example.} Let us consider the system
$$
\mat{c}{\dot x\hl z\\ e}=
\mat{c|cc}{A&B&B_d\hl C&D&D_{zd}\\ C_e&0&0}
\mat{c}{x\hl w\\d},\
\mat{c|cc}{A&B&B_d\hl C&D&D_{zd}\\ C_e&0&0}:=
{\small\left(\begin{array}{cccc|cc}
-0.97&2.2&2.36&3.45&-0.62&-0.1\\
-0.21&-0.8&5.2&-0.35&-0.7&-0.32\\
-2.56&-4.97&-0.75&-9.75&-1.42&-0.84\\
-3.64&0.2&9.68&-0.64&0&0\\ \hline
0&-0.36&0.36&-0.57&-1.14&-1.76\\ \hdashline
1.5&-0.11&0&0.93&0&0\\
0.1&0&0&0&0&0
\end{array}\right)}
$$
in feedback with $w=\Delta_\delta z$ where $\delta\in[-0.6,5]$ and $x(0)=0$ for simplicity. For computing a ``smallest'' ellipsoid that contains the output trajectory $e(.)$ against disturbances $d(.)$ whose energy is bounded by one, we  minimize the trace of $Y$ over the LMIs
\eqn{
\kyp{
\Xf
}{
\mat{ccccc}{
\mi   &0\\
0     &-I}}{
\mat{cc|cccc}{
A_\Psi&B_\Psi T\Cgr&
B_\Psi T\Dgr&
B_\Psi T\Ddp \\0&A&B&B_d\hl
C_\Psi&D_\Psi T\Cgr& D_\Psi T\Dgr&D_\Psi T \Ddp\\
0&0&0&I}
}\cl 0\te{with}\Ddp:=\mat{c}{D_{zd}\\0},
}
\eqn{
\mat{c|ccc}{
Y&0&0&C_e\hl
0&\Xf_{11}&\Xf_{12}-K&\Xf_{13}\\
0&\Xf_{21}-K^T&\Xf_{22}&\Xf_{23}\\
C_e^T&\Xf_{31}&\Xf_{32}&\Xf_{33}}\cg 0\te{and}
\mat{ccc|ccc|cc}{
R_{11} &R_{12}-K     \\
R_{21}-K^T &R_{22}   }\cl 0;
}
here we work with $\psi_1$, $\psi_2$ as in Remark~\ref{rem2} for $\nu=0,1,2,3$.
If the LMIs are feasible, straightforward adaptations of the proof of Theorem \ref{Tmain} reveal that all trajectories of the uncertain system satisfy $e(T)^TY^{-1}e(T)\leq \int_0^T d(t)^Td(t)\,dt$ for $T>0$; if $\|d\|\leq 1$, this translates into the ellipsoidal invariance property $e(T)\in\c E:=\{e\in\R^2\st e^TY^{-1}e\leq 1\}$ for $T>0$. The numerical results are depicted in Fig.~\ref{fig1}. The blue ellipsoid is obtained for static multipliers ($\nu=0$), while the red, yellow and black ones ($\nu=1,2,3$) clearly exhibit the benefit of the dynamics in the multipliers. The one for $\nu=3$ is
 tight, as supported by
interconnection trajectories for the worst parameter value $\delta=-0.6$ and five worst-case disturbance
inputs (with energy one) hitting the boundary of the black ellipsoid at different points.


\ls{
\mun{
\bul
\mat{cc|ccc}{
0&0&\Xf_{11}&\Xf_{12}\\
0&0&\Xf_{21}&\Xf_{22}\hl
\Xf_{11}&\Xf_{12}&0&0\\
\Xf_{21}&\Xf_{22}&0&0}
\mat{cc|cccc}{
I&0&0&0\\
0&I&0&0\hl
A_\Psi&B_\Psi T\Cgr&
B_\Psi T\Dgr&
B_\Psi TD_d\\0&A&B&B_d
}+\\+\bul
\mat{ccccc}{
\mi&0\\
0  &-I}
\mat{cc|cccc}{
C_\Psi&D_\Psi T\Cgr& D_\Psi T\Dgr&D_\Psi D_d\\
0&0&0&I}
\cl 0
}
\mun{
\mat{c|cc}{A&B&B_d\hl C&D&D_d\\0&I&0\\C_e&0&0}
\mat{ccc}{I&0\\\delta(I-D\delta)^{-1}C&\delta(I-D\delta)^{-1}D_d\\0&I}=\\=
\mat{c|cc}{
A+B\delta(I-D\delta)^{-1}C&B_d+B\delta(I-D\delta)^{-1}D_d\hl C+D\delta(I-D\delta)^{-1}C &D_d+\delta(I-D\delta)^{-1}D_d\\
\delta(I-D\delta)^{-1}C&\delta(I-D\delta)^{-1}D_d\\C_e&0}=\\=
\mat{c|cc}{A+B\delta(I-D\delta)^{-1}C&B_d+B\delta(I-D\delta)^{-1}D_d\hl (I-D\delta)^{-1}C &(I-D\delta)^{-1}D_d\\\delta(I-D\delta)^{-1}C&\delta(I-D\delta)^{-1}D_d\\C_e&0}=
\mat{c|cc}{A_\delta&B_\delta\hl \mat{c}{1\\\delta}C_\delta & \mat{c}{1\\\delta}D_\delta\\C_e&0}
}
\mun{
\bul
\mat{cc|ccc}{
0&0&\Xf_{11}&\Xf_{12}\\
0&0&\Xf_{21}&\Xf_{22}\hl
\Xf_{11}&\Xf_{12}&0&0\\
\Xf_{21}&\Xf_{22}&0&0}
\mat{cc|cccc}{
I&0&0\\
0&I&0\hl
A_\Psi&B_\Psi T\mat{c}{1\\\delta}C_\delta&
B_\Psi T\mat{c}{1\\\delta}D_\delta\\0&A_\delta&B_\delta
}+\\+\bul
\mat{ccccc}{
\mi&0\\
0  &-I}
\mat{cc|cccc}{
C_\Psi&D_\Psi T \mat{c}{1\\\delta}C_\delta&D_\Psi T\mat{c}{1\\\delta}D_\delta\\
0&0&I}
\cl 0
}
$$
\arr{rcl}{\dot x&=&Ax+Bw+B_dd,\ \ x(0)=0,\\z&=&C x+D w+D_dd\\z&=&C x+D w+D_dd\\w&=&0x+Iw+0d\\e&=&C_ex}
$$
$$
\arr{rcl}{
\mat{c}{\dot \xi\\\dot x}&=&
\mat{cc}{
A_\Psi&B_\Psi\Cgr\\
0     &A         }
\mat{c}{\xi\\x}+
\mat{cc}{
B_\Psi\Dgr& B_\Psi D_d\\
B&0}\mat{c}{w\\d},\ \ \mat{c}{\xi(0)\\x(0)}=0\\
\mat{c}{y\\z\\d}&=&
\mat{cc}{
C_\Psi&D_\Psi\Cgr\\
0&\Cg\\
0&0}\mat{c}{\xi\\x}+
\mat{cccccc}{
D_\Psi\Dgr&D_\Psi D_d\\
\Dg&I_\nz\\
0&I_\nz}\mat{c}{w\\d};
}
$$
$$
\mat{cccc}{
\Xf_{11}&\Xf_{12}-K&\Xf_{13}\\
\Xf_{21}-K^T&\Xf_{22}&\Xf_{23}\\
\Xf_{31}&\Xf_{32}&\Xf_{33}}
\cg
\mat{cccc}{0\\0\\C_e^T}H^{-1}
\mat{cccc}{0&0&C_e}.
$$
$$
x(T)^TC_e^TH^{-1}C_ex(T)\leq
\mat{c}{\xi_1(T)\\\xi_2(T)\\x(T)}^T
\mat{cccc}{
\Xf_{11}&\Xf_{12}-K&\Xf_{13}\\
\Xf_{21}-K^T&\Xf_{22}&\Xf_{23}\\
\Xf_{31}&\Xf_{32}&\Xf_{33}}
\mat{c}{\xi_1(T)\\\xi_2(T)\\x(T)}\leq
\int_0^T\!  \|d(t)\|^2 \, dt
$$
$$
\mas{c|c}{A_\delta&B_\delta\hl C&0}
$$
$$
A_\delta W+WA_\delta^T+B_\delta B_\delta^T=0
$$
States reachable with energy one:
$$
\{x\st x^TW^{-1}x=1\}
$$
Choose $H:=CWC^T$ to infer
$$
\mat{cc}{H&C\\C^T&W^{-1}}\cge 0\te{and thus}W^{-1}\cge C^TH^{-1}C.
$$
Kernel vectors $x_0$ of $W^{-1}\cge C^TH^{-1}C$ with $x_0^TW^{-1}x_0=1$ lead to outputs $z_0=Cx_0$
with $z_0^TH^{-1}z_0=1$.

$$
\dot x(t)=(A_\delta-B_\delta F)x(t),\ z(t)=Cx(t),\ u(t)=-Fx(t),\ x(0)=\xi,\ t\in[0,T]
$$
$$
\dot x(t)=A_\delta x(t)+B_\delta u(t),\ z(t)=Cx(t),\ x(0)=\xi_0,\ x(T)=\xi_T,\ t\in[0,T]
$$
For $\t x(t)=x(T-t)$ and $\t u(t)=u(T-t)$ we infer
$$
\dot {\t x}(t)=-\dot x(T-t)=-A_\delta x(T-t)-B_\delta u(T-t),\ z(t)=Cx(t),\ x(0)=\xi,\ t\in[0,T]
$$
}{}

\section{Conclusions}
In this paper, we have given a complete link between the general IQC theorem and dissipation theory. To this end we proposed a new stability result based on the notion of finite horizon IQCs with a terminal cost. For a classical frequency domain test related to
parametric uncertainties, it {was}
shown that one can work with convexly constrained terminal cost matrices,
the benefit of which {was}
illustrated through a numerical example.
It is hope{d} that this framework lays the foundation for further research on guaranteeing local time-domain properties through tests emerging in absolute stability theory.

%
%
\section*{References}
\bibliographystyle{elsart-num}
{\footnotesize
\bibliography{ref}
}

\appendix

\section{An auxiliary result}\label{apndx1}

\lemma{\label{Lau1}For compatibly sized matrices (with $T,R,S$ being invertible) suppose that
$$
\mat{cc}{\t A&\t B\\\t C&\t D}=
\mat{cc}{T^{-1}&0\\0&R^{-1}}\mat{cc}{A&B\\C&D}\mat{cc}{T&0\\F&S}.
$$
Then
$$\arc{0.3}
\kyp{X}{\mi}{\mat{cc}{A&B\\C&D}}\cl 0
\ \Leftrightarrow\
\kyp{T^TXT}{R^T\mi R}{\mat{cc}{\t A&\t B\\\t C&\t D}}\cl 0.
$$

}
{\bf Proof.} The equivalence follows since the left-hand sides of the inequalities are related by a congruence transformation with $\diag(T,S)$.\epro

\section{Proof of Lemma \ref{Lfdi2}}\label{apndx2}
Given $\t Z$ and any certificate $\Xf$ of \r{fdi}, we need to prove $\Xf-\diag(\t Z,0)\cg 0$.
Let us choose a certificate $Y=Y^T$ of \r{lup}.
\ls{
$$
T^{-T}\mat{cc}{-\t\mi&0\\0&\t\mi}T^{-1}=
\frac{1}{4}\mat{cc}{-I& I\\I& I}
\mat{cc}{\t\mi &-\t\mi \\\t\mi &\t\mi}=
\frac{1}{2}\mat{cc}{0&\t\mi\\\t\mi&0},
$$
$$
T\mat{cc}{\t\Psi E&0\\ 0&\t\Psi F}=
\mat{cc}{\t\Psi\mat{cc}{-E&F}\\ \t\Psi\mat{cc}{E&F}}=
\mat{cc}{\t \Psi\t G\t \Psi^{-1}\\ I}\t \Psi\mat{cc}{E&F},
$$
and thus
\eqn{
\mat{cc}{-\t\mi&0\\0&\t\mi}=T^T\mat{cc}{0&\frac{1}{2}\t\mi\\\frac{1}{2}\t\mi&0}T\te{and}
T\mat{cc}{\t\Psi E&0\\ 0&\t\Psi F}=
\mat{cc}{\t \Psi\t G\t \Psi^{-1}\\ I}\t \Psi\mat{cc}{E&F}.}
}
{}
With
$\mat{cc}{\t \Cf&\t \Df}:=\mat{cc|c}{\t C_\Psi&\t D_\Psi\Cgr&\t D_\Psi\Dgr}
$,
the KYP inequalities corresponding to \r{fdi} and \r{lup} for $\Pi=\t\Psi^*\t \mi\t\Psi$ read as
\eql{kyph}{
\mat{cccccc}{
I&0\\
\Af&\Bf\\
\t\Cf&\t\De}^T
\mat{ccc}{
0&\t\Xf&0\\
\t\Xf&0&0\\
0&0&\t\mi}\!
\mat{cccccc}{
I&0\\
\Af&\Bf\\
\t\Cf&\t\De}\cl 0\te{and}
\mat{cccccc}{
\bul\\
\bul\\
\bul}^T
\mat{ccc}{
0&\t\Y&0\\
\t\Y&0&0\\
0&0&\t\mi}\!
\mat{cccccc}{
I&0\\
A_\Psi&B_{\Psi}E\\
\t C_\Psi&\t D_{\Psi}E}\cg 0.
}
Right-multiplying \r{are1} with $\diag\left(I,\mat{cc}{\Cgr&\Dgr}\right)$ and left-multiplying the transpose shows
\eql{are2}{
\mat{cccccc}{
I&0\\
\Af&\Bf\\
\Cf&\De}^T
\mat{ccc}{
0&\Zf&0\\
\Zf&0&0\\
0&0&\mi}\!
\mat{cccccc}{
I&0\\
\Af&\Bf\\
\Cf&\De}=
\mat{ccc}{\t \Cf&\t \Df}^T\t\mi\mat{ccc}{\t \Cf&\t \Df}
}
for $\Zf:=\diag(\t Z,0_n)$.
If we subtract \r{are2} from \r{kyp}, we infer that $\t\Xf:=\Xf-\Zf$ satisfies the first LMI
in \r{kyph}.
For $Y$ certifying \r{lup}, one argues analogously to see that $\t Y:=Y-\t Z$
satisfies the second LMI in \r{kyph}. These two inequalities can be diagonally combined to
$$
\kyp{
\mat{cccccc}{
-\t Y       &0             &0            \\
0           &\t \Xf_{11}   &\t \Xf_{12}  \\
0           &\t \Xf_{21}   &\t \Xf_{22}  }
}
{
\mat{ccc|ccc}{
-\t\mi &0    \\
0      &\t\mi}}
{\mat{ccc|ccc}{
A_\Psi      &0           &0                &B_{\Psi}E     &0               \\
0           &A_\Psi      &B_\Psi \Cgr      &0             &B_\Psi\Dgr      \\
0           &0           &A                &0             &B               \hl
\t C_{\Psi} &0           &0                &\t D_{\Psi}E  &0               \\
0           &\t C_{\Psi} &\t D_{\Psi}\Cgr  &0             &\t D_{\Psi}\Dgr }}\cl 0.
$$
Note that this certifies \r{fdic}. Let us now see how the move from \r{fdic} to \r{fdit} proceeds for this LMI based on Lemma \ref{Lau1}. We start with easily verified equation
\mun{\arc{0.3}
\mat{ccc|cc}{-I&I&0&0&0\\I&I&0&0&0\\0&0&I&0&0\hl 0&0&0&-I&I\\0&0&0&I&I}
\mat{ccc|ccc}{
A_\Psi      &0           &0                &B_{\Psi}E    &0               \\
0           &A_\Psi      &B_\Psi \Cgr      &0            &B_\Psi\Dgr      \\
0           &0           &A                &0            &B               \hl
\t C_{\Psi} &0           &0                &\t D_{\Psi}E &0               \\
0           &\t C_{\Psi} &\t D_{\Psi}\Cgr  &0            &\t D_{\Psi}\Dgr }
\mat{ccc|cc}{
-\frac{1}{2}I&\frac{1}{2}I&0&0&0\\\frac{1}{2}I&\frac{1}{2}I&0&0&0\\0&0&I&0&0\hl 0&0&-C&I&-D\\0&0&0&0&I}
\mat{ccc|cc}{I&0&0&0\\0&I&0&0\\0&0&I&0\hl -\t D_\Psi^{-1}\t C_\Psi&0&0&\t D_\Psi^{-1}}
=
\\=
\mat{ccc|ccc}{
A_\Psi      &-B_\Psi\mat{cc}{-E&\Dgr}\t D_{\Psi}^{-1}\t C_\Psi   &2B_\Psi \Cgr     &B_\Psi\mat{cc}{-E&\Dgr}\t D_{\Psi}^{-1}  \\
0           &A_\Psi-B_\Psi\t D_\Psi^{-1}\t C_\Psi                &0                & B_\Psi\t D_{\Psi}^{-1}                      \\
0           &-B_\Psi\mat{cc}{0 &B}\t D_{\Psi}^{-1}\t C_\Psi      &A                & \mat{cc}{0&B}\t D_{\Psi}^{-1}               \hl
\t C_{\Psi} &-D_\Psi\mat{cc}{-E&\Dgr}\t D_{\Psi}^{-1}C_\Psi  &2\t D_{\Psi}\Cgr &D_\Psi\mat{cc}{-E&\Dgr}\t D_{\Psi}^{-1}  \\
0           &0           &0                & I                           }.}
(The result forms indeed is a realization of $\col(H,I)$, but this is not relevant for the arguments that follow.) In view of the first equation in \r{tr1} and by Lemma \ref{Lau1}, we conclude
$$
\kyp{K}{\mat{cc}{0&\frac{1}{2}\t\mi\\\frac{1}{2}\t\mi&0}
}{\mat{c|c}{A_H&B_H\hl C_H&D_H\\0&I}}\cl 0
$$
for
\eqn{
K:=\mat{ccccc}{-\frac{1}{2}I&\frac{1}{2}I&0\\\frac{1}{2}I&\frac{1}{2}I&0\\0&0&I}^T
\mat{cccc}{
-\t Y &0         &0      \\
0     &\t \Xf_{11} &\t \Xf_{12} \\
0     &\t \Xf_{21} &\t \Xf_{22} }
\mat{ccccc}{-\frac{1}{2}I&\frac{1}{2}I&0\\\frac{1}{2}I&\frac{1}{2}I&0\\0&0&I}=
\mat{ccccc}{
\frac{1}{4}(\t \Xf_{11}-\t Y)&\frac{1}{4}(\t \Xf_{11}+\t Y) &\frac{1}{2}\t \Xf_{12}\\
\frac{1}{4}(\t \Xf_{11}+\t Y)&\frac{1}{4}(\t \Xf_{11}-\t Y) &\frac{1}{2}\t \Xf_{12}\\
\frac{1}{2}\t \Xf_{21}       &\frac{1}{2}\t \Xf_{21}        &\t \Xf_{22}}.
}
The left-upper block of the LMI reads as $A_H^TK+KA_H\cl 0$. By inspection, $A_H$ is Hurwitz.
Therefore, $K\cg 0$, which in turn shows $\t\Xf=\Xf-\Zf\cg 0$ (and also $\t Y\cl 0$, i.e., $\t Z\cg Y$).

\section{Proof of Lemma \ref{Liqc}}\label{apndx3}
Choose any $z\in\Le_{2e}^\nz$ and $T>0$ and define
$$
\t y:=\t\Psi\mat{c}{z\\\Delta(z)}\in\Le_{2e}^{\nz+\nw}\te{as well as}
\mat{c}{\t z\\\t w}:=\t\Psi^{-1}(\t y_T)\in\Le_2^{\nz+\nw}.
$$
By causality of $\t\Psi^{-1}$ we infer
$$
\mat{c}{\t z_T\\\t w_T}=
\mat{c}{\t z\\\t w}_T=(\t\Psi^{-1}(\t y_T))_T =(\t\Psi^{-1} \t y)_T= \mat{c}{z\\\Delta(z)}_T=\mat{c}{z_T\\\Delta(z)_T}.
$$
Hence $\col(\t z,\t w)$ is identical to $\col(z,\Delta(z))$ on $[0,T]$ and constitutes a modification of the latter trajectory on $(T,\infty)$ in order to generate a finite energy signal. As the crucial point,
this modification even has the property
\eql{in1}{
\int_0^T \left[\t\Psi\mat{c}{z\\\Delta(z)}\right]^T\!\!\!\t\mi\t\Psi\mat{c}{z\\\Delta(z)}\,dt=
\int_0^\infty \left[\t\Psi\mat{c}{\t z\\\t w}\right]^T\!\!\!\t\mi\t
\Psi\mat{c}{\t z\\\t w}\,dt\geq 0
}
and can hence be interpreted as a stable extension in the sense of Yakubovich \cite{Yak00,Yak02}.

To prove \r{in1} we first observe through a simple computation (using bilinearity) that
\mun{
\arc{.3}
\int_0^\infty \left[\t\Psi\mat{c}{\t z\\2\Delta(\t z)-\t w}\right]^T\!\!\!\t\mi\t\Psi\mat{c}{\t z\\\t w}\,dt=\\=
\arc{.3}
\int_0^\infty \left[\t\Psi\mat{c}{\t z\\\Delta(\t z)}\right]^T\!\!\!\t\mi\t\Psi\mat{c}{\t z\\\Delta(\t z)}\,dt
-
\int_0^\infty \left[\t\Psi\mat{c}{0\\\Delta(\t z)-\t w}\right]^T\!\!\!\t\mi\t\Psi\mat{c}{0\\\Delta(\t z)-\t w}\,dt.}
If we exploit \r{iqc} and \r{rl} we thus conclude
\eql{in2}{\arc{.3}
\int_0^\infty \left[\t\Psi\mat{c}{\t z\\2\Delta(\t z)-\t w}\right]^T\!\!\!\t\mi\t\Psi\mat{c}{\t z\\\t w}\,dt\geq 0.
}
By causality of $\Delta$ we have $\t w_T=\Delta(z)_T=\Delta(z_T)_T=\Delta(\t z_T)_T=\Delta(\t z)_T$, and
causality of $\t\Psi$ then implies
\eqn{
\left[\t\Psi\mat{c}{\t z\\2\Delta(\t z)-\t w}\right]_T=
\left[\t\Psi\mat{c}{\t z_T\\2\Delta(\t z)_T-\t w_T}\right]_T=
\left[\t\Psi\mat{c}{\t z_T\\\t w_T}\right]_T=
\left[\t\Psi\mat{c}{z_T\\\Delta(z)_T}\right]_T.
}
This leads to
\eql{eq3}{
\left[\t\Psi\mat{c}{\t z\\2\Delta(\t z)-\t w}\right]_T=
\left[\t\Psi\mat{c}{z\\\Delta(z)}\right]_T=\t y_T=\t\Psi\mat{c}{\t z\\\t w}.
}
Since the signal \r{eq3} is supported on $[0,T]$, we arrive at
\eqn{
\int_0^\infty \left[\t\Psi\mat{c}{\t z\\2\Delta(\t z)-\t w}\right]^T\!\!\!\t\mi\t\Psi\mat{c}{\t z\\\t w}\,dt=
\int_0^\infty \left[\t\Psi\mat{c}{\t z\\2\Delta(\t z)-\t w}\right]_T^T\!\!\!\t\mi \t y_T\,dt=
\int_0^T \t y_T^T\t\mi \t y_T\,dt.}
Due to \r{eq3}, this equals the left-hand side of \r{in1}, which in turn proves
\r{in1} by \r{in2}.

In a final step we consider
$
y:=\Psi\mat{c}{z\\\Delta(z)}\in\Le_{2e}^\ny.
$
Since $y$ is the response of \r{psi} for $u=\col(z,\Delta(z))$ and since
$\t y=\t C_\Psi\xi+\t D_\Psi u$, we can right-multiply \r{are1} with $\col(\xi,u)$ and left-multiply the transpose to infer
$
\frac{d}{dt} \xi(t)^T\t Z\xi(t)+y(t)^T\mi y(t)=\t y(t)^T\t\mi\t y(t)\te{for}t\geq 0.
$
By integration and using $\xi(0)=0$ we get
$$
\xi(T)^T\t Z\xi(T)+\int_0^T y(t)^T\mi y(t)\,dt=\int_0^T\t y(t)^T\t\mi\t y(t)\,dt\te{for}T\geq 0.
$$
The combination with \r{eq3} and \r{in1} concludes the proof.

\section{Proof of Lemma \ref{Liqc2}}\label{apndx4}

{\bf Proof.} By \r{grl}, the new multiplier
\eql{defm}{
\Pi_H:=
\mat{cc}{I&H\\0&I}^*
\Pi
\mat{cc}{I&H\\0&I}=\Psi_H^*\mi\Psi_H\te{with}\Psi_H:=\Psi\mat{cc}{I&H\\0&I} 
}
satisfies \r{rl}. With a realization $(A_H,B_H,C_H,D_H)$ of $H$ where $A_H$ is Hurwitz,
we get
$$
\mat{cc}{I&H\\0&I}=\mas{c|cc}{A_H&0&B_H\hl C_H&I&D_H\\0&0&I}=:\mas{c|cc}{A_H&B_e\hl C_e&D_e}
\te{and}
\Psi\mat{cc}{I&H\\0&I}=
\mas{cc|cccc}{
A_\Psi&B_{\Psi}C_e &B_{\Psi}D_e\\
0     &A_H         &B_e        \hl
C_\Psi&D_{\Psi}C_e &D_{\Psi}D_e}.
$$
Note that $A_H-B_eD_e^{-1}C_e=A_H$ is Hurwitz as well.
If we right-multiply \r{are1}
with the matrix $\diag(I,(C_e\ D_e))$ and left-multiply with the transpose, we obtain
\eql{are3}{
\mat{ccc}{\bul\\\bul\\\bul\\\bul\\\bul}^T
\mat{cc|cc|c}{0&0&\t Z&0&0\\0&0&0&0&0\hl \t Z&0&0&0&0\\0&0&0&0&0\hl 0&0&0&0&\mi }
\mat{cc|c}{I&0&0\\0&I&0\hl A_\Psi&B_\Psi C_e&B_\Psi D_e\\0     &A_H         &B_e        \hl C_\Psi&D_\Psi C_e&D_\Psi C_e}=
(\bul)^T \mi \mat{cc|c}{\t C_\Psi&\t D_\Psi C_e&\t D_\Psi D_e}.
}
As argued in Section \ref{Scf}, we infer that $\diag(\t Z,0)$ hence certifies the  factorization
$
\Pi_H=\t \Psi_H^*\t\mi\t \Psi_H\te{with}
\t \Psi_H=
\mas{cc|cccc}{
A_\Psi    & B_{\Psi}C_e &B_{\Psi}D_e\\
0         &A_H          &B_e        \hl
\t C_\Psi &\t D_\Psi C_e&\t D_\Psi D_e }.
$
This is even a canonical one since both the state-matrix of $\t\Psi_H$ and the following one of the inverse $\t\Psi_H^{-1}$ are Hurwitz, as seen by inspection:
$$
\mat{cccccc}{
A_\Psi    & B_{\Psi}C_e \\
0         &A_H          }-
\mat{cccccc}{
B_{\Psi}D_e\\
B_e        }(\t D_\Psi D_e)^{-1}
\mat{cccccc}{
\t C_\Psi & \t D_\Psi C_e }=
\mat{cccccc}{
A_\Psi-B_{\Psi}\t D_\Psi^{-1}\t C_\Psi&\bul             \\
0                       &A_H-B_eD_e^{-1}C_e }.
$$

Next choose any $z\in\Le_2^\nz$ and define
\eql{eq4}{
\mat{c}{v\\w}:=\mat{cc}{I&-H\\0&I}\mat{cc}{z\\\Delta(z)}=
\mat{cc}{(I-H\Delta)(z)\\\Delta(z)}.
}
With the causal and bounded system $\Delta_H:=\Delta(I-H\Delta)^{-1}$ we infer from
$v=(I-H\Delta)(z)$ that $z=(I-H\Delta)^{-1}(v)$ and hence
$w=\Delta(z)=\Delta_H(v)$, thus implying
\eql{tra}{
\mat{cc}{I&H\\0&I}\mat{c}{v\\\Delta_H(v)}=\mat{cc}{I&H\\0&I}\mat{c}{v\\w}=\mat{cc}{z\\\Delta(z)}.
}
Due to \r{defm} and by \r{iqc1} we hence get
\eqn{
\int_{-\infty}^\infty
\mat{c}{\pl\\\pl }^*\Pi_H(\io)
\mat{c}{\h v(\io)\\ \widehat{\Delta_H(v)}(\io) }\,d\om=
\int_{-\infty}^\infty
\mat{c}{\pl\\ \pl }^*\Pi(\io)
\mat{c}{\hat z(\io)\\ \widehat{\Delta(z)}(\io) }\,d\om\geq 0.
}

All this allows us to apply Lemma~\ref{Liqc} for $\Delta_H$ and the multiplier $\Pi_H$ with the canonical factorization \r{defm} as certified by $\diag(\t Z,0)$. Hence, for any $v\in\Le_{2e}^\nz$, the
response of $y_H:=\Psi_H\col(v,\Delta_H(v))$ with state-trajectory $\col(\xi,x_H)$
satisfies
\eql{dis2}{
\int_0^T y_H(t)^T\mi y_H(t)\,dt+\mat{c}{\xi(T)\\x_H(T)}^T\mat{cc}{\t Z&0\\0&0}\mat{c}{\xi(T)\\x_H(T)}
\geq 0\te{for}T>0.
}

Finally, let $y$ be the output of \r{psi} driven by $u=\col(z,\Delta(z))$ with $z\in\Le_{2e}^\nz$.
Then \r{eq4} defines $\col(v,w)\in\Le_{2e}^{\nz+\nw}$ and \r{tra} persists to hold. This implies
$$
y_H=\Psi_H\mat{c}{v\\\Delta_H(v)}=
\Psi\mat{c}{z\\\Delta(z)}=y,
$$
by which \r{dis2} is identical to \r{disu} for $Z:=\t Z$ and which completes the proof.

\section{Proof of Lemma \ref{LK}}\label{apndx5}

Feasibility of \lmiRa\ with $M\in\c M$ implies that $\psi_1^*P\psi_2$ is GPR. Due to \cite{Hel95,Goh96b} this guarantees the existence of a (non-symmetric) canonical factorization of $\psi_1^*P\psi_2$ certified by the solution $K$ of
$A_1^TK+KA_2+C_1^TPC_2-(KB_2+C_1^TPD_2)(D_1^TPD_2)^{-1}(B_1^TK+D_1^TPC_2)=0$
under the stability constraints
$\eig(A_1^T-(KB_2+C_1^TPD_2)(D_1^TPD_2)^{-1}B_1^T)\subset\Cl$
and $\eig(A_2-B_2(D_1^TPD_2)^{-1}(B_1^TK+D_1^TPC_2))\subset\Cl$.
One checks by a direct computation that $\t Z$ defined by \r{cer}
is indeed a certificate for $\Psi^*M\Psi$ having a canonical factorization
$\t\Psi^*\t M\t\Psi$; since $T$ is invertible, the same holds for
$\Pi=[\Psi T]^*M[\Psi T]$ and $[\t\Psi T]^*\t M[\t\Psi T]$.

The statement on $\Delta\in\bm{\Delta}$
the follows from Lemma \ref{Liqc}. Since $\Pi$ is a PN mutliplier, Lemma \ref{Lfdi2} guarantees  \coub. The same line of reasoning applied to \lmiRa\ (with $G=I$ for $-M$, $-R$, $-K$) leads to \coub.

\end{document}